\documentclass[11pt,twoside]{article}
\usepackage{times}
\usepackage{amsmath,amssymb,amsthm}
\usepackage{enumerate}
\usepackage{cite}
\usepackage{graphicx}
\usepackage{epstopdf}
\usepackage{color}
\usepackage{cases}

\allowdisplaybreaks

\def \no{\nonumber}

\newcommand{\pf}{\noindent {\bf Proof. \hspace{2mm}}}
\newcommand{\ef}{ \hfill $ \square $ \vskip 3mm}
\newcommand{\be}{\begin{equation}}
\newcommand{\ee}{\end{equation}}
\newcommand{\bea}{\begin{eqnarray}}
\newcommand{\eea}{\end{eqnarray}}
\def\p{\partial}
\def\ve{\varepsilon}

\def\f{\frac}
\def\na{\nabla}
\def\la{\lambda}
\def\al{\alpha}

\def\vp{\varphi}

\def\si{\sigma}

\def\dl{\delta}

\def\ds{\displaystyle}

\def\dl{\delta}
\def\no{\nonumber}

\def\beq{\begin{equation}}
\def\eeq{\end{equation}}
\def\ben{\begin{eqnarray}}
\def\een{\end{eqnarray}}
\def\bec{\begin{cases}}
\def\eec{\end{cases}}
\def\ss{\sum\limits}

\begin{document}
\footskip=0pt
\footnotesep=2pt
\let\oldsection\section
\renewcommand\section{\setcounter{equation}{0}\oldsection}
\renewcommand\thesection{\arabic{section}}
\renewcommand\theequation{\thesection.\arabic{equation}}
\newtheorem{claim}{\noindent Claim}[section]
\newtheorem{theorem}{\noindent Theorem}[section]
\newtheorem{lemma}{\noindent Lemma}[section]
\newtheorem{proposition}{\noindent Proposition}[section]
\newtheorem{definition}{\noindent Definition}[section]
\newtheorem{remark}{\noindent Remark}[section]
\newtheorem{corollary}{\noindent Corollary}[section]
\newtheorem{example}{\noindent Example}[section]
\title{Formation and construction of large variational shock waves for 1-D $n\times n$ quasilinear hyperbolic
conservation systems}
\author{Yin Huicheng \quad Zhu Wanqing\footnote{Yin Huicheng
(huicheng@nju.edu.cn, 05407@njnu.edu.cn) and Zhu Wanqing (zwq15371032536@163.com) are supported by the NSFC
(No.~12331007) and by the National key research and development program of China (No. 2020YFA0713803).}\vspace{0.5cm}\\
\small
School of Mathematical Sciences and Institute of Mathematical Sciences, \\
\small Nanjing Normal University, Nanjing 210023, P. R. China.
}
\vspace{0.5cm}

\date{}
\maketitle

\centerline {\bf Abstract} \vskip 0.3 true cm
In the paper [Li Jun, Xu Gang, Yin Huicheng, On the blowup mechanism of smooth solutions to 1D quasilinear strictly
hyperbolic systems with large variational initial data, Nonlinearity 38 (2025), No.2, 025016],
for the 1-D $n\times n$ ($n\geqslant 3$) strictly hyperbolic system  $\p_tv+F(v)\p_xv=0$ with some classes of
large variational initial data $v(x, 0)$,
the geometric blowup mechanism and the detailed singularity behaviours
of  $\p_{x,t}v$ near the blowup point are studied when the $n\times n$ matrix $F(v)$ admits at least one genuinely nonlinear
eigenvalue. In this paper, we focus on the formation and construction of a large variational shock wave from the blowup point for
1-D $n\times n$ quasilinear hyperbolic conservation law system  $\p_tu+\p_xf(u)=0$ when some smooth simple wave solution
is generic non-degenerate before the formation of singularity and the corresponding eigenvalue is genuinely nonlinear.

\vskip 0.3 true cm
{\bf Keywords:} Hyperbolic conservation law, genuinely nonlinear, shock,  entropy condition,

\qquad\qquad\quad generic non-degenerate, Rankine-Hugoniot condition.

\vskip 0.3 true cm
\tableofcontents

\section{Introduction}
Consider the Cauchy problem of 1-D $n\times n$ quasilinear strictly hyperbolic system:
\begin{subequations}\label{zhu}
\begin{numcases}{}
\p_tu+\p_xf(u)=0,\label{1.1-1}\\[2mm]
u(x,-1)=u_0(x),\label{1.1-2}
\end{numcases}
\end{subequations}
where $n\geqslant 2$, $t\geqslant -1$, $x\in\Bbb R$, $u=(u_1, \cdots, u_n)^{\top}: \mathbb{R}\times[-1,\infty)\to \mathbb{R}^n$ is the unknown function
vector of $(x,t)$ and $f=(f_1, \cdots, f_n)^{\top}: \mathbb{R}^n\to \mathbb{R}^n$ is smooth in its arguments. The initial data $u_0(x)
=(u_{1,0}(x), \cdots, u_{n,0}(x))\in C_0^4(\mathbb{R})$ and $\|u_0(x)\|_{L^{\infty}(\mathbb{R})}$ is small but $\|u_0^{(m)}(x)\|_{L^{\infty}(\mathbb{R})}$
($m=1,\cdots, 4$) may be large.
By the strict hyperbolicity of \eqref{1.1-1}, the $n\times n$ matrix $F(u)=\nabla_u f(u)$ admits $n$ distinct real eigenvalues
\begin{equation}\label{utezheng}
  \la_1(u)<\cdots<\la_n(u),
\end{equation}
whose right eigenvectors are denoted by $r_1(u), \cdots, r_n(u)$, respectively.
It is assumed that  \eqref{1.1-1} is either linearly degenerate or genuinely nonlinear with respect to
each eigenvalue $\la_k(u)$ ($1\leqslant k\leqslant n$), namely,
$\na_u\la_k(u)\cdot r_k(u)\equiv 0$ or $\na_u\la_k(u)\cdot r_k(u)\neq0$ holds.
In addition, it follows from Proposition 2.1 of \cite{LXY} that
when \eqref{1.1-1} is genuinely nonlinear with respect to some fixed $\la_i(u)$ ($1\leqslant i\leqslant n$),
there exists an invertible mapping $u\mapsto w=(w_1, \cdots, w_n)^{\top}=(g_1(u), \cdots, g_n(u))^{\top}=g(u)$ with  $g(0)=0$ such that for $C^1$ smooth solution
$u$ ($|u|$ is small but $|\p_{x,t}u|$ may be large), \eqref{zhu} can be equivalent to
\begin{subequations}\label{zhuyao}
\begin{numcases}{}
\p_t w+A(w)\p_x w=0,\label{1.3-1}\\[2mm]
w(x,-1)=w_0(x)=(w_{1,0}(x),\cdots, w_{n,0}(x)),\label{1.3-2}
\end{numcases}
\end{subequations}
where $A(w)=(\p_ug)F(g^{-1}(w))(\p_ug)^{-1}:=(a_{ij}(w))_{n\times n}$. Meanwhile, the following properties hold
\begin{itemize}
  \item $A(w)$ has $n$ distinct eigenvalues: $\la_1(w), \cdots, \la_n(w)$, where $\la_k(w)=\la_k(g^{-1}(w))$ for $1\leqslant k\leqslant n$;
  \item the elements of $A(w)$ satisfy $a_{ji}(w)=0$ for $j\neq i$ and $a_{ii}(w)=\la_i(w)$;
  \item $\p_{w_{i}}\la_i(0)\neq 0$;
  \item $A(0)=diag(\la_1(0),\cdots,\la_n(0))$;
   \item $w_0(x)\in C_0^4(\Bbb R)$.
\end{itemize}

Let $l_j(w)$ be the left eigenvector of $A(w)$ corresponding to the eigenvalue $\la_j(w)(1\leqslant j\leqslant n)$.
Without loss of generality, one can assume that
\begin{align*}
& l_j(w)\cdot r_k(w)=\delta^k_j,\quad 1\leqslant j,k\leqslant n; \\
& r_i(w)=\vec{e_i},\quad\text{$r_j(0)=\vec{e_j}$ and $\|r_j(w)\|=1$ for  $j\neq i$}; \\
& l_j(0)=\vec{e_j}^\top \quad\text{for $1\leqslant j\leqslant n$},
\end{align*}
where $\|r_j(w)\|=\sqrt{\sum\limits^n_{k=1}(r_{j,k}(w))^2}$ with $r_j(w)=(r_{j,1}(w),\cdots, r_{j,n}(w))^\top$.

We now study the $i$-simple wave solution $\overline{w}=(0,\cdots,0, w_i, 0, \cdots, 0)$ of \eqref{1.3-1}
with the initial data $\overline{w}_0(x)=(0,\cdots,0, w_{i,0}(x), 0, \cdots, 0)\in C_0^4(\Bbb R)$
\beq\label{isimple}
\p_t w_i+\la_i(\overline{w})\p_x w_i=0,
\eeq
where
\begin{equation}\label{supp}
w_i(x,-1)=w_{i,0}(x)\not\equiv 0,\quad\text{and $|w_{i,0}(x)|$ is small but $|w'_{i,0}(x)|$ may be large}.
\end{equation}
Set
\begin{equation}\label{1-1}
G(x):=\la_i(\overline{w}_0(x)),
\end{equation}
then $G(x)-G(0)\in C_0^4(\mathbb{R})\bigcap L^{\infty}(\mathbb{R})$.
Note that the $C^4$ smooth solution $w_i$ of \eqref{isimple} will blow up at the maximal existence time $T_*=-1-\dfrac{1}{\min_{x\in \mathbb{R}} G'(x)}$.
We assume that $x=0$ is the unique and generic non-degenerate minimal point of $G'(x)$ for $x\in\mathbb{R}$.
Set $G(0)=\la_i(0)=0$ (otherwise, the translation $(x,t)=(x+\la_i(0)t,t)$ can be done). Without loss of generality, suppose that
\begin{equation}\label{1-2}
G'(0)=\min_{x\in \mathbb{R}}G'(x)=-1,\quad G''(0)=0,\quad G'''(0)=6.
\end{equation}
In addition, $\p_{w_i}\la_i(0)=1$ is assumed (otherwise, one can take the spatial scaling $(x,t)=(\p_{w_i}\la_i(0)x,t)$).
At this time, it holds that
\begin{align}\label{itz}
\la_i(\overline{w})=w_i+O(w_i^2),
\end{align}
$T_*=0$ is the blowup time of smooth solution $w_i$   to \eqref{isimple} with \eqref{supp}
and \eqref{1-2}, and the corresponding unique blowup point is $(x_*,T_*)=(0,0)$.

Our main result is:
\begin{theorem}\label{main1}
Let $u_0(x)=g^{-1}(\overline{w}_0(x))$ in \eqref{1.1-2} and \eqref{1-2} hold. Assume that
\eqref{1.1-1} is genuinely nonlinear with respect to some fixed $\la_i(u) (1\leqslant i\leqslant n)$.
Then there exists a small positive constant $\ve>0$ such that problem \eqref{zhu}  admits a
weak entropy solution $u\in C^{\f13}(\Bbb R\times [0,\ve])$ with an $i$-shock curve $x=\varphi(t)\in C^2[0, \ve]$
starting from the blowup point $(0,0)$, and $(\vp(t), u(x,t))$
fulfills the following estimations for $t\geqslant 0$:
\begin{align}\label{1.12-0}
&\varphi(t)=\la_i(u(0,0))t+O(t^2),\no\\
&u(x,t)=u(0,0)+O((t^3+(x-\la_i(u(0,0))t)^2)^{\f{1}{6}}).
\end{align}
Meanwhile, corresponding to the solution $w\in C^{\f13}(\Bbb R\times [0,\ve])$ of \eqref{zhuyao}, we have that for $t\geqslant 0$,
\begin{align}\label{1.12}
&w_i(x,t)=O((t^3+(x-\varphi(t))^2)^{\f{1}{6}}),\no\\
&w_j(x,t)=O(t(t^3+(x-\varphi(t))^2)^{\f{1}{6}}),\quad\quad j\neq i.
\end{align}
\end{theorem}

\begin{remark}\label{CYY-1}
Choosing $w_{i,0}(x)=\delta^{\nu}\phi(\f{x}{\dl})$ with $\dl>0$ being small and
$\phi(s)\in C_0^{4}(\Bbb R)$ but $\phi(s)\not\equiv 0$, where $0<\nu<1$ is some fixed constant,
then for the corresponding initial data $u_0(x)$ of \eqref{zhu}, it holds that $|u_0(x)|$ is small but  $|u'_0(x)|$ is large.
In addition, for the initial data $\overline{w}_0(x)=(C_1,\cdots, C_{i-1}, C_i+w_{i,0}(x), C_{i+1}, \cdots, C_n)\in C^4(\Bbb R)$
with $(C_1,\cdots, C_n)$ being some constant vector, we have the analogous conclusions as in Theorem \ref{main1}.
\end{remark}

\begin{remark}\label{CYY-2}
When $w_0(x)=\overline{w}_0(x)+\dl {\widetilde w}_0(x)$ is the small perturbation of $\overline{w}_0(x)$,
where $\dl>0$ is small and ${\widetilde w}_0(x)\in C_0^{\infty}(\Bbb R)$, it follows from the result in \cite{A0}
that the geometric blowup of \eqref{zhuyao} is stable under the generic non-degenerate condition \eqref{1-2}.
In this case, by collecting the methods in this paper and \cite{DY}, the analogous conclusions in Theorem \ref{main1}
still hold.
\end{remark}

\begin{remark}\label{CYY-2-0}
If \eqref{zhu} can be transformed into \eqref{zhuyao} for large $|u|$
(for example, this is valid for the 1-D $3\times 3$ full compressible Euler equations of polytropic gases with no
vacuum), Theorem \ref{main1} still holds even for large initial data $|u_0(x)|$.
However, for large $|u|$, one can not obtain the form \eqref{zhuyao} from \eqref{zhu} in general.
\end{remark}

\begin{remark}\label{CYY-3}
Theorem \ref{main1} can be directly applied to many physical models such as 1-D compressible Euler equations of polytropic gases
($3\times3$ system), 2-D supersonic steady full compressible Euler equations ($4\times4$ system), 1-D MHD equations ($5\times5$ system),
1-D elastic wave equations ($6\times6$ system) and 1-D full ideal compressible MHD equations ($7\times7$ system).
One can see these related models in \cite{LXY}, \cite{YSL} and the references therein.
For example, note that the equations of 3-D ideal compressible magnetohydrodynamics (MHD) are
\begin{equation}\label{FFF-6}
\begin{cases}
\ds\p_t\rho+div(\rho u)=0,\\
\ds\p_t(\rho u)+div(\rho u\otimes u-H\otimes H)+\nabla(P+\f12|H|^2)=0,\\
\ds\p_tH-curl(u\times H)=0,\\
\ds div H=0,\\
\ds \p_t(\rho S)+div(\rho uS)=0,\\
\end{cases}
\end{equation}
where $(x,t)=(x_1,x_2,x_3,t)$, $\rho$ is the fluid density,
$u=(u_1,u_2,u_3)^{\top}$ is the fluid velocity, $H=(H_1,H_2,H_3)^{\top}$
is the magnetic field, $S$ is the entropy
and $P$ is the pressure satisfying the state equation $P=P(\rho,S)=A\rho^{\gamma}e^{\f{S}{c_v}}$
with $A, c_v$ and $\gamma>1$ being positive constants.
Let $(\rho, u, H, S)(x,t)=(\rho, u, H, S)(x_1,t)$ and $H_1=\bar H_1>0$
is a constant. Then \eqref{FFF-6} becomes the 1-D conservation law
\begin{equation}\label{FFF-6-1}
\begin{cases}
\ds\p_t\rho+\p_1(\rho u_1)=0,\\
\ds\p_t(\rho u_1)+\p_1(\rho u_1^2)+\p_1(P+\f12|H_2|^2+\f12|H_3|^2)=0,\\
\ds\p_t(\rho u_2)+\p_1(\rho u_1 u_2-\bar H_1H_2)=0,\\
\ds\p_t(\rho u_3)+\p_1(\rho u_1 u_3-\bar H_1H_3)=0,\\
\ds\p_tH_2+\p_1(u_1H_2-\bar H_1u_2)=0,\\
\ds\p_tH_3+\p_1(u_1H_3-\bar H_1u_3)=0,\\
\ds \p_t(\rho S)+\p_1(\rho u_1S)=0.\\
\end{cases}
\end{equation}
For the $C^1$ smooth solution $(\rho, u_1, u_2, u_3, H_2, H_3, S)$,  \eqref{FFF-6-1} is equivalent to
\begin{equation}\label{FFF-6-2}
\begin{cases}
\ds\partial_t\rho+u_1\p_1\rho+\rho\p_1u_1=0,\\[5pt]
\ds\partial_tu_1+u_1\p_1u_1+\f{c^2}{\rho}\p_1\rho+\f{\p_SP}{\rho}\p_1S+\f{H_2}{\rho}\p_1H_2+\f{H_3}{\rho}\p_1H_3=0,\\[5pt]
\ds\partial_tu_2+u_1\p_1u_2-\f{\bar H_1}{\rho}\p_1H_2=0,\\[5pt]
\ds\partial_tu_3+u_1\p_1u_3-\f{\bar H_1}{\rho}\p_1H_3=0,\\[5pt]
\ds\partial_tH_2+u_1\p_1H_2-\bar H_1\p_1u_2+H_2\p_1u_1=0,\\[5pt]
\ds\partial_tH_3+u_1\p_1H_3-\bar H_1\p_1u_3+H_3\p_1u_1=0,\\[5pt]
\p_tS+u_1\p_1S=0.
\end{cases}
\end{equation}
The initial data of \eqref{FFF-6} is given by
\begin{equation}\label{FFF-7}
\begin{split}
&\rho(x, 0)=\bar\rho+\rho_0(x),\quad u_1(x, 0)=u_{1,0}(x), \quad u_2(x, 0)=u_{2,0}(x),\quad u_3(x, 0)=u_{3,0}(x),\\
&H_2(x, 0)=\bar H_2+H_{2,0}(x), \quad H_3(x, 0)=\bar H_3+H_{3,0}(x), \quad S(x, 0)=\bar S+S_0(x),
\end{split}
\end{equation}
where $\bar\rho, \bar H_2, \bar H_3$ and $\bar S$ are constants, $\bar\rho>0$, $\bar H_2\bar H_3\not=0$,
$(\rho_0(x), u_{1,0}(x), u_{2,0}(x), u_{3,0}(x), H_{2,0}(x),
$ $H_{3,0}(x),S_0(x))\in C_0^{4}(\Bbb R)$,  $|\rho_0(x)|+\sum\limits_{j=1}^3|u_{j,0}(x)|+\sum\limits_{j=2}^3|H_{j,0}(x)|
+|S_0(x)|$ is small but $|\rho'_0(x)|+\sum\limits_{j=1}^3|u'_{j,0}(x)|$ $+\sum\limits_{j=2}^3|H'_{j,0}(x)|
+|S'_0(x)|$ may be large.

For $(\rho, u_1, u_2, u_3, H_2, H_3, S)=(\bar\rho,0,0,0,\bar H_2,\bar H_3,\bar S)$,
\eqref{FFF-6-2} has seven real distinct eigenvalues
\begin{equation*}
\begin{split}
&\ds\lambda_1=-\bigg\{\f{|\bar H|^2}{2\bar\rho}+\f{{\bar c}^2}{2}
+\f12\sqrt{(\f{|\bar H|^2}{\bar\rho}+{\bar c}^2)^2-\f{4}{\bar\rho}{\bar H}_1^2{\bar c}^2}
\bigg\}^{\f12}\\
&<\lambda_2=-\sqrt{\f{1}{\bar\rho}}{\bar H}_1\\
&<\lambda_3=-\bigg\{\f{|\bar H|^2}{2\bar\rho}+\f{{\bar c}^2}{2}
-\f12\sqrt{(\f{|\bar H|^2}{\bar\rho}+{\bar c}^2)^2-\f{4}{\bar\rho}{\bar H}_1^2{\bar c}^2}
\bigg\}^{\f12}\\
&<\lambda_4=0\\
&<\lambda_5=-\lambda_3<\lambda_6=-\lambda_2<\lambda_7=-\lambda_1,
\end{split}
\end{equation*}
where $|\bar H|^2={\bar H}_1^2+{\bar H}_2^2+{\bar H}_3^2$ and $\bar c=\sqrt{\p_{\rho}P(\bar\rho,\bar S)}$.
It can be verified directly that \eqref{FFF-6-2} is genuinely nonlinear with respect to all the eigenvalues except $\lambda_4$
for small perturbations of $(0,0,0,\bar\rho,\bar H_2,\bar H_3,\bar S)$.
Thus, for the $i$-simple wave solution of \eqref{FFF-6} ($i\not=4$),
Theorem \ref{main1} is suitable for the formation and construction of $i$-shock.
\end{remark}

\begin{remark}\label{CYY-4}
There have been a series  of interesting results  on the shock formation of smooth solutions and the descriptions
of pre-shocks to the multidimensional compressible Euler equations or the second order nonlinear wave
equations of divergence type (see \cite{A1}-\cite{A0}, \cite{Buck1}-\cite{Buck4}, \cite{Ch1}, \cite{Ch2},
\cite{Ding-1}-\cite{Ding-2},\cite{0-Speck},
\cite{Lu1}-\cite{MY} and \cite{Sp}), which illustrate that the
formation of the multidimensional shock is due to the compression of
the characteristic surfaces.
In addition,  for the multidimensional compressible perfect Euler equations,
the author in \cite{Ch4} constructs a restricted shock (but the corresponding Rankine-Hugoniot  conditions
are not satisfied). That is, the related construction of a multidimensional
shock wave after the blowup of smooth solution is not obtained yet except the scalar multidimensional
quasilinear hyperbolic conservation law (see \cite{YZ2}).
\end{remark}

We now recall some results on the formation and construction of shocks for 1-D hyperbolic conservation laws.
For 1-D scalar convex conservation law,  \cite{CZ} proved the  formation  and construction
of shocks for piecewise smooth initial data with finite discontinuities. For the general 1-D scalar conservation law, \cite{YZ1}
described the shock formation and development under the various assumptions of initial data.
For 1-D $2\times2$ $p$-system of polytropic gases, under the assumptions that one Riemann invariant is a constant and the initial data satisfy the related generic
non-degenerate condition for another Riemann invariant,
\cite{Lebaud} constructed a shock solution from the blowup point. This result was extended to the more
general case of $p$-system in \cite{ChenD, Kong}, where both the  Riemann invariants are not constants and
it is additionally supposed that only one family of characteristics is  squeezed, while the other characteristics family does not
squeeze at the same point.
For 1-D $3\times3$ strictly hyperbolic conservation law with small generic non-degenerate initial data,
\cite{CXY} constructed the weak entropy shock solution from the blowup point. Recently,
\cite{DY} improved the result in \cite{CXY} to the more general 1-D $n\times n(n\geqslant3)$ strictly or symmetric
hyperbolic conservation laws with small generic non-degenerate initial data.
For the 3-D full compressible Euler equations with spherically symmetric solutions of small perturbations
of non-vacuum constant states, \cite{Y} and \cite{Ch3} investigate the shock formation and development.
For some classes of large variational initial data, the authors in \cite{LXY}
demonstrated that for the 1-D $n\times n$ strictly hyperbolic system, the geometric blowup
can happen in finite time meanwhile the solution itself still remains bounded, especially,
for the 1-D $n\times n$ strictly hyperbolic conservation law, the shock is formed.

In this paper, we mainly focus on the construction of a shock for 1-D $n\times n$ strictly hyperbolic conservation law
\eqref{zhu} with the $i$-simple wave initial data. At first, at the blowup time $T_*=0$,
we give the precise fractional series expansion on the pre-shock $i$-simple wave of \eqref{isimple},
and singular $w_i(x,0)$ will be taken as the initial data for the shock construction of problem \eqref{zhuyao}.
Secondly, in order to avoid the involved analysis on the blowup system in \cite{CXY} and \cite{DY} for small
data solution, motivated by the methods in \cite{Buck1} for 2-D compressible Euler equations with azimuthal symmetry,
we will solve the singular initial data problem and singular initial-boundary value problem of \eqref{zhuyao} on two sides
of some fixed approximate shock curve $x=\phi(t)$ starting from the point $(0,0)$ and give the detailed time-space estimates
on the singularities of $w$ and its derivatives.
In this process, the good form of $A(w)$ in \eqref{zhuyao},
the Lax geometric entropy conditions, the analysis on related different characteristics
and the choice of suitable iteration schemes play basic roles.
Based on this, the real shock curve $x=\vp(t)$ with $\vp(0)=0$
can be found by the contraction mapping principle and subsequently Theorem \ref{main1} is shown.

Our paper are organized as follows: In Section \ref{Sec-2}, the shock formation for the $i$-simple wave is investigated
and the related  fractional expansions of the pre-shock with the $C^{\f{1}{3}}$ H$\mathrm{\ddot{o}}$lder regularity near
the blowup point are derived in details.
In Section \ref{Sec-3}, in terms of the entropy conditions,
the suitable iteration schemes to solve \eqref{zhuyao} are chosen by the corresponding singular initial value problem and singular initial-boundary
value problems, and the resulting solution $w$ of \eqref{zhuyao} can be obtained on two sides of the fixed approximate shock curve.
In Section \ref{Sec-4}, the $C^2$ regular shock curve corresponding to the solution of \eqref{zhuyao} is constructed,
moreover, we establish the uniqueness of the weak entropy solution in $C^{\f13}(\Bbb R\times [0,\ve])$
with some additional singularities of derivatives near the blowup point $(0,0)$.
Based on Section \ref{Sec-2}-Section \ref{Sec-4}, the main results in Theorem \ref{main1}
are derived in Section \ref{Sec-5}.

\section{Properties of the $i$-simple wave}\label{Sec-2}
In this section, for reader's convenience and applications later,
we will give the details for the shock formation and development of the $i$-simple wave \eqref{isimple}
and establish the related precise estimates of pre-shock although some properties have been obtained
in \cite{YZ1} and the references therein.

The characteristics $x=\eta(\beta, t)$ of \eqref{isimple} starting from $(\beta, -1)$ is defined as
\begin{equation*}
\begin{cases}
\p_t\eta(\beta,t)=\la_i(\overline{w}(\eta(\beta,t),t)), \\[2mm]
\eta(\beta,-1)=\beta.
\end{cases}
\end{equation*}
It follows from \eqref{isimple} and \eqref{1-2} that
\begin{equation}\label{et}
\eta(\beta,t)=\beta+(t+1)\la_i(\overline{w}_0(\beta)),~~w_i(\eta(\beta, t), t)=w_0(\beta),
\end{equation}
the blowup time and the blowup point are $T_*=0$ and  $(x_*, T_*)=(0,0)$, respectively.

By \eqref{1-1}-\eqref{1-2} and \eqref{et}, it holds that
\begin{equation*}
\eta(\beta,t)\in C^4(\mathbb{R}\times[-1,0])
\end{equation*}
and
\begin{equation*}
\p_{\beta}\eta(0,0)=0,\quad\p^2_{\beta}\eta(0,0)=0,\quad\p^2_{t\beta}\eta(0,0)=-1,\quad\p^3_{\beta}\eta(0,0)=6.
\end{equation*}

Then near $\beta=0$, one has
\begin{align}\label{tlzkphi}
\eta(\beta,0)&=\eta(0,0)+\p_{\beta}\eta(0,0)\beta+\f{1}{2}\p^2_{\beta}\eta(0,0)\beta^2
+\f{1}{6}\p^3_{\beta}\eta(0,0)\beta^3+\f{1}{24}\p^4_{\beta}\eta(\beta_1,0)\beta^4\no\\
& =\beta^3+\f{1}{24}\p^4_{\beta}\eta(\beta_1,0)\beta^4,\quad \text{where $\beta_1$ lies in $0$ and $\beta$.}
\end{align}
Setting $\eta(\beta,0)=x$ and using Lemma 4.14 of \cite{Buck1}, we obtain that from \eqref{tlzkphi}
\begin{equation}\label{xtheta}
  \beta=x^{\f{1}{3}}+\al_2x^{\f{2}{3}}+\al_3x+O(|x|^{\f{4}{3}}),
\end{equation}
where
\begin{equation*}
  \al_2=-\f{1}{3}\cdot\f{1}{24}\p^4_{\beta}\eta(\beta_1,0), ~~\al_3=\f{1}{3}\Big(\f{1}{24}\p^4_{\beta}\eta(\beta_1,0)\Big)^2.
\end{equation*}
In addition,  it is derived from \eqref{et}, \eqref{itz} and \eqref{1-2} that near $\beta=0$,
\begin{align}\label{w0}
w_i(\eta(\beta,0),0)&=w_0(\beta)\no\\
& =-\beta+\f{1}{2}w''_0(0)\beta^2+\f{1}{6}w'''_0(0)\beta^3+\f{1}{24}w^{(4)}_0(\beta_1)\beta^4.
\end{align}
Combining \eqref{xtheta} and \eqref{w0} yields
\begin{align}
  w_i(x,0) & =-x^{\f{1}{3}}+(-\al_2+\f{1}{2}w''_0(0))x^{\f{2}{3}}\no\\
   & ~~~~+(-\al_3+\al_2w''_0(0)+\f{1}{6}w'''_0(0))x+O(|x|^{\f{4}{3}})\no\\
   &:=-x^{\f{1}{3}}+a_2x^{\f{2}{3}}+a_3x+O(|x|^{\f{4}{3}}).\label{w013}
\end{align}

Next, we estimate $\p_xw_i(\eta(\beta,0),0)$.
Taking the first order derivative with respect to $\beta$ on both sides of the second term in \eqref{et} and using
the Taylor expansion of  $ \p_{\beta}\eta(\beta,0)$, one has
\begin{equation}\label{2-8}
 \begin{cases}
\p_xw_i(\eta(\beta,0),0)\p_{\beta}\eta(\beta,0)=w'_0(\beta)=-1+w''_0(0)\beta+O(\beta^2), \\[2mm]
    \p_{\beta}\eta(\beta,0)=3\beta^2+\f{1}{6}\p^4_{\beta}\eta(\beta_2,0)\beta^3, \quad \text{where $\beta_2$ lies in $0$ and $\beta$.}
 \end{cases}
\end{equation}
This, together with \eqref{xtheta}, yields
\begin{align}\label{w1}
  \hspace{-0.4cm}\p_xw_i(\eta(\beta,0),0) &=\f{w'_0(\beta)}{\p_{\beta}\eta(\beta,0)}\no\\
  &=-\f{1}{3}x^{-\f{2}{3}}+\Big(\f{2}{3}\al_2+\f{1}{3}w''_0(0)+\f{1}{54}\p^4_{\beta}\eta(\beta_2,0)\Big)x^{-\f{1}{3}}+O(1).
\end{align}
By comparing the coefficients of $\p_xw_i(\eta(\beta,0),0)$ in \eqref{w1} and the expression
of $w_i'(x,0)$ from \eqref{w013}, we arrive at
\begin{equation*}
\f{2}{3}\al_2+\f{1}{3}w''_0(0)+\f{1}{54}\p^4_{\beta}\eta(\beta_2,0)=\f{2}{3}a_2.
\end{equation*}
Hence,
\begin{equation*}
  \p_xw_i(x,0)=-\f{1}{3}x^{-\f{2}{3}}+\f{2}{3}a_2x^{-\f{1}{3}}+O(1).
\end{equation*}
In order to treat $\p^2_xw_i(\eta(\beta,0),0)$, differentiating the first expression in \eqref{2-8} with respect to $\beta$
and using the Taylor expansion of  $ \p_{\beta}^2\eta(\beta,0)$  give
\begin{equation}\label{p2xw}
\begin{cases}
\p^2_xw_i(\eta(\beta,0),0)(\p_{\beta}\eta(\beta,0))^2+\p_xw_i(\eta(\beta,0),0)\p^2_{\beta}\eta(\beta,0)=w''_0(0)+O(\beta),\\[2mm]
\p^2_{\beta}\eta(\beta,0)=6\beta+\f{1}{2}\p^4_{\beta}\eta(\beta_3,0)\beta^2, \quad \text{where $\beta_3$ lies in $0$ and $\beta$.}
\end{cases}
\end{equation}
Together  with \eqref{w1} and \eqref{p2xw}, by direct computation, one has
\begin{align*}
  \p^2_xw_i(\eta(\beta,0),0)& =\f{w''_0(0)-\p_xw_i(\eta(\beta,0),0)\p^2_{\beta}\eta(\beta,0)}{(\p_{\beta}\eta(\beta,0))^2} \no\\
   & =\f{2}{9}\beta^{-5}+O(\beta^{-4})=\f{2}{9}x^{-\f{5}{3}}+O(x^{-\f{4}{3}}).
\end{align*}
Analogously,
\begin{align*}
  &\p^3_xw_i(\eta(\beta,0),0)\no\\
  =&\f{w'''_0(\beta)-\p_xw_i(\eta(\beta,0),0)\p^3_{\beta}\eta(\beta,0)
  -3\p^2_xw_i(\eta(\beta,0),0)\p_{\beta}\eta(\beta,0)\p^2_{\beta}\eta(\beta,0)}{(\p_{\beta}\eta(\beta,0))^3}\no\\
  =&O(\beta^{-8})=O(x^{-\f{8}{3}}).
\end{align*}
Near $x=O(\ve)$, it holds that
\begin{align}
  &|w_i(x,0)|\leqslant M ,\label{cz0}\\
  &|w_i(x,0)+x^{\f{1}{3}}-a_2x^{\f{2}{3}}|\leqslant M|x|,\label{cz01}\\
  &|w'_i(x,0)+\f{1}{3}x^{-\f{2}{3}}-\f{2}{3}a_2x^{-\f{1}{3}}|\leqslant M,\label{cz1}\\
  &|w''_i(x,0)-\f{2}{9}x^{-\f{5}{3}}|\leqslant M|x|^{-\f{4}{3}},\label{cz2}\\
  &|w'''_i(x,0)|\leqslant M|x|^{-\f{8}{3}},\label{cz3}
\end{align}
where the constant $|a_2|\leqslant M$ and $M\gg1$.

It is pointed out that in the process of constructing a shock starting from the blowup point $(0,0)$,
the properties \eqref{cz0}-\eqref{cz3} will play important roles.
Assume that the shock curve is denoted by $x=\varphi(t)$ with $\varphi(0)=0$.
Then it follows from the Rankine-Hugoniot conditions of \eqref{1.1-1} that
 \begin{equation}\label{RHu}
\begin{cases}
\si(t)[u_i(w)]=[f_i(u(w))],\\[2mm]
\si(t)[u_j(w)]=[f_j(u(w))], \qquad j\neq i,
\end{cases}
\end{equation}
where $\si(t)=\varphi'(t)$, $[u_i(w)]=u_i(w_-)-u_i(w_+)$, and the $\mp$   stand for
the states before or behind the shock curve. The corresponding Lax geometric entropy conditions are given by (see
Chapter 15 of \cite{Smoller})
\begin{equation}\label{eq:3.28-1}
\lambda_{i-1}(u(w_{-}))<\sigma<\lambda_i(u(w_{-})), \quad \lambda_i(u(w_+))<\sigma<\lambda_{i+1}(u(w_+)).
\end{equation}
Note that by \eqref{RHu}, as in (4.25) of \cite{DY}, we have
\begin{equation}\label{RH-0}
\si=\la_i\Big(\int^1_0\p_{u_k}f_l(\theta u(w_+)+(1-\theta)u(w_-))d\theta\Big)
\end{equation}
with $\la_i(\int^1_0\p_{u_k}f_l(\theta u(w_+)+(1-\theta)u(w_-))d\theta)$ being
the $i$-th eigenvalue of $n\times n$ matrix \\
$\Big(\int^1_0(\partial_{u_k}f_l)(\theta u(w_+)+
 (1-\theta)u(w_-))d\theta\Big)_{k,l=1}^n$
and \begin{equation}\label{RH}
  [w_j]=\mathcal{F}_j(w_{1,+},\cdots,w_{i-1,+}, w_{i,\pm}, w_{i+1,-},\cdots,w_{n,-})[w_i]^3,\quad j\neq i,
\end{equation}
where $\mathcal{F}_j$ are smooth functions on their arguments.

Firstly, it is assumed that the approximate shock curve $x=\phi(t)\in C^2([0,\ve])$ satisfies
\begin{equation}\label{jb}
   |\phi(t)|\leqslant M^4t^2,~~~~~~~~~~|\phi'(t)|\leqslant 2M^4t,~~~~~~~~~~|\phi''(t)|\leqslant M^9,
\end{equation}
where the constant $M>0$ is sufficiently large. Then given the following initial data for \eqref{1.3-1}
\begin{equation}\label{2-19}
    (w_{1,0}(x),\cdots,w_{n,0}(x))=(0,\cdots,0,w_i(x,0),0,\cdots,0),
\end{equation}
here $w_i(x,0)$ satisfies \eqref{cz0}-\eqref{cz3},
the approximate shock solution of \eqref{zhuyao} can be constructed for $x=\phi(t)$ in \eqref{jb} with the jump discontinuity \eqref{RH}, but $x=\phi(t)$ is not the solution of \eqref{RH-0} in general. Actually, the real shock curve $x=\varphi(t)$ solves \eqref{RH-0} and satisfies \eqref{jb}, which will be determined later.

Next, under the condition \eqref{jb}, we construct the pre-shock solution of \eqref{zhuyao} from \eqref{cz0}-\eqref{cz3}.
Let
 \begin{equation}\label{w20}
\begin{cases}
\p_tw^{(0)}_i+\la_i(\overline{w}^{(0)})\p_xw^{(0)}_i=0,~~~~~~~~~in~~~~\mathcal{B}_{\ve}=([-M\ve,M\ve]\times[0,\ve])\backslash\{(\phi(t),t):
t\in[0,\ve]\},\\[2mm]
 w^{(0)}_i=w_i(x,0),\quad\quad\quad\quad\quad\quad\quad~~for~~~~x\in [-M\ve,M\ve].
\end{cases}
\end{equation}
Then it follows from the characteristics method that for $t\in[0,\ve]$,
\beq\label{eta0}
\text{$w^{(0)}_i(\eta^{(0)}(\beta,t),t)=w_i(\beta,0)$ and $\eta^{(0)}(\beta,t)=\beta+t\la_i(\overline{w}(\beta,0))$.}
\eeq

\begin{lemma}\label{0}
Under the condition of \eqref{jb}, there exists
a largest $\beta^{(0)}_+(t)>0$ and a smallest $\beta^{(0)}_-(t)<0$ for \eqref{eta0} such that
\beq\label{v}
\phi(t)=\eta^{(0)}(\beta^{(0)}_{\pm}(t),t),\quad\quad\quad\f{4}{5}t^{\f{3}{2}}<|\beta^{(0)}_{\pm}(t)|<\f{6}{5}t^{\f{3}{2}}.
\eeq
\end{lemma}
\pf  Set
\beq\label{g0}
g_0(\beta)=\la_i(\overline{w}(\beta,0))+\beta^{\f{1}{3}}-\Big(\f{1}{2}\p^2_{w_i}\la_i(0)+a_2\Big)\beta^{\f{2}{3}}.
\eeq
By \eqref{itz}, \eqref{cz01} and \eqref{cz1}, we have $|g_0(\beta)|\leqslant2M^2|\beta|$ and $|g'_0(\beta)|\leqslant2M^2$.
For $t>0$, let
\beq\label{t}
\tau=t^{\f{1}{2}},\quad\quad y=\beta^{\f{1}{3}}\tau^{-1},\quad\quad\zeta=\phi(t)\tau^{-3}.
\eeq
In terms of \eqref{jb}, $|\zeta|\leqslant M^4\tau\ll1$ holds.
Note that the first equation in \eqref{v} comes from the equality
\begin{align}\label{z}
\zeta=&y^3-y+\Big(\f{1}{2}\p^2_{w_i}\la_i(0)+a_2\Big) y^2\tau+\tau^{-1}g_0(y^3\tau^3)\no\\
:=&y^3-y+\mathcal{L}(y,\tau).
\end{align}
Define $c=\f{1}{2}\p^2_{w_i}\la_i(0)+a_2$ for convenience. Then for $|y|\leqslant10$ and $0<\tau\leqslant\ve$,
one has
\begin{subequations}
\begin{align}
&|\mathcal{L}(y,\tau)-c y^2\tau|\leqslant 2M^3\tau^2,\label{L-1}\\
&|\p_y\mathcal{L}(y,\tau)-2cy\tau|\leqslant 6M^3\tau^2,\label{L-2}\\
&|\p_{\tau}\mathcal{L}(y,\tau)-cy^2|\leqslant 8M^3\tau.\label{L-3}
\end{align}
\end{subequations}

By \eqref{z} and \eqref{L-1},
when $\tau=0$, we have $\zeta=y^3-y$ and
\begin{equation*}
y_{\pm}(\zeta,0)=\pm1+\f{1}{2}\zeta\mp\f{3}{8}\zeta^2+\f{1}{2}\zeta^3\mp\f{105}{128}\zeta^4+\f{3}{2}\zeta^5+O(|\zeta|^6)\quad
\text{for~~$|\zeta|\leqslant\f{1}{2}$.}
\end{equation*}
Taking the first order derivative with respect to $y$ on both sides of \eqref{z}  yields for $\tau=0$
\begin{align*}
\p_y(y^3-y+\mathcal{L}(y,\tau)-\zeta)|_{(y_{\pm}(\zeta,0),0)}=3y^2_{\pm}(\zeta,0)-1\neq0.
\end{align*}
By the implicit function theorem, $y_{\pm}=y_{\pm}(\zeta,\tau)$ is obtained from \eqref{z}. Moreover,
for $|\zeta|\leqslant M^4\tau\ll1$, it follows from \eqref{z}, \eqref{L-2} and \eqref{L-3} that
\begin{align*}
y_{\pm}(\zeta,\tau)&=y_{\pm}(\zeta,0)+\p_{\tau} y_{\pm}(\zeta,0)\tau+O(\tau^2)=\pm1+\f{1}{2}\zeta-\f{c}{2}\tau+O(\tau^2).
\end{align*}
From the definitions in \eqref{t}, we set
\begin{equation}\label{2-26}
\beta^{(0)}_{\pm}(t)=y_{\pm}^3(\zeta,\tau)\tau^3=(\pm1+\f{1}{2}\zeta-\f{c}{2}\tau+O(\tau^2))^3\tau^3,
\end{equation}
which implies
\begin{equation}\label{2-30}
|\beta^{(0)}_{\pm}(t)\mp t^{\f{3}{2}}|\leqslant M^4t^2.
\end{equation}
Then \eqref{v} is shown.
\ef

Set $\Upsilon_0(t)=[-M\ve,M\ve]\backslash[\beta^{(0)}_-(t),\beta^{(0)}_+(t)]$.
By \eqref{cz0}-\eqref{cz3}, it holds that for $\beta\in\Upsilon_0(t)$,
\begin{align}\label{czt}
   & |w_i(\beta,0)|\leqslant M\ve^{\f{1}{3}},~\quad\quad\quad\quad~~-\f{2}{5}t^{-1}\leqslant w'_i(\beta,0)\leqslant -\f{1}{5}t^{-1},\no\\
   &|w''_i(\beta,0)|\leqslant\f{1}{3}t^{-\f{5}{2}},\quad\quad\quad~~~~~~|w'''_i(\beta,0)|\leqslant 2Mt^{-4}.
\end{align}
Furthermore, it follows from $\f{9}{10}\leqslant\p_{w_i}\la_i(\overline{w}(\beta,0))\leqslant\f{11}{10}$ in \eqref{itz} and \eqref{czt} that
\beq\label{fm}
1+t\p_{\beta}\la_i(\overline{w}(\beta,0))\geqslant\f{1}{2},\quad\quad
\Big|\f{w'_i(\beta,0)}{1+t\p_{\beta}\la_i(\overline{w}(\beta,0))}\Big|\leqslant2t^{-1}.
\eeq
Define the map
\begin{equation*}
\eta^{(0)}(\cdot,t):\Upsilon_0(t)\longrightarrow[-M\ve,M\ve]\backslash\{\phi(t)\},
\end{equation*}
we then have
\begin{lemma}\label{eta}
It holds that for $s\in[0,t]$ and $\beta\in\Upsilon_0(t)$,
\begin{align}
&|\p_{\beta}\eta^{(0)}(\beta,s)-1|\leqslant\f{3}{5},\label{b1}\\
&|\eta^{(0)}(\beta,s)-\phi(s)|\geqslant\f{3}{5}t^{\f{1}{2}}(t-s),\label{b2}
\end{align}

and the inverse map $\eta^{(0),-1}(\cdot,t):[-M\ve,M\ve]\backslash\{\phi(t)\}\longrightarrow\Upsilon_0(t)$
fulfills
\begin{align}
&\f{5}{8}\leqslant\p_x\eta^{(0),-1}(x,t)\leqslant\f{5}{2},\label{a}\\
&\f{4}{5}t^{\f{3}{2}}+\f{5}{8}|x-\phi(t)|\leqslant|\eta^{(0),-1}(x,t)|\leqslant\f{6}{5}t^{\f{3}{2}}+\f{5}{2}|x-\phi(t)|.\label{c}
\end{align}
\end{lemma}
\pf
Taking the derivative with respect to $\beta$ on both sides of \eqref{eta0} yields
\begin{equation*}
\p_{\beta}\eta^{(0)}(\beta,s)=1+s\p_{\beta}\la_i(\overline{w}(\beta,0)).
\end{equation*}
Together with \eqref{itz} and \eqref{czt}, one has $|\p_{\beta}\eta^{(0)}(\beta,s)-1|\leqslant\f{3}{5}$.
Then $\eta^{(0)}$ is an injection and is a strictly increasing function on $\beta\in\Upsilon_0$.

On the other hand, by
$\phi(t)=\eta^{(0)}(\beta^{(0)}_{-}(t),t)=\eta^{(0)}(\beta^{(0)}_{+}(t),t)$, \eqref{cz01}, \eqref{jb}, \eqref{v} and \eqref{b1},
we have $\eta^{(0)}(\beta^{(0)}_{+}(t),s)>\phi(s)$ for $s\in(0,t)$ and the following lower bound holds
\begin{align}\label{phixj}
  \eta^{(0)}(\beta^{(0)}_{+}(t),s)-\phi(s) & =-\int^t_s(\p_{t'}\eta^{(0)}(\beta^{(0)}_{+}(t),t')-\phi'(t'))dt' \no\\
   &=\int^t_s(\phi'(t')-\la_i(\overline{w}(\beta^{(0)}_{+}(t),0)))dt'\no\\
   &\geqslant(\beta^{(0)}_{+}(t)^{\f{1}{3}}-2|c|t)(t-s)-M^4(t^2-s^2)\no\\
   & \geqslant\f{3}{5}t^{\f{1}{2}}(t-s).
\end{align}
Similarly,  $\phi(s)-\eta^{(0)}(\beta^{(0)}_{-}(t),s)\geqslant\f{3}{5}t^{\f{1}{2}}(t-s)$ can be achieved. Thus \eqref{b2} is shown.

In addition, for $(x,t)\in\mathcal{B}_{\ve}$, due to $\p_x\eta^{(0),-1}
=\f{1}{\p_{\beta}\eta^{(0)}}$
and
\begin{equation*}
\eta^{(0),-1}(x,t)-\beta^{(0)}_{\pm}(t)=\eta^{(0),-1}(x,t)-\eta^{(0),-1}(\phi(t),t)=(x-\phi(t))(\p_x\eta^{(0),-1})(\bar{x},t),
\end{equation*}
where $\bar{x}$ lies in $x$ and $\phi(t)$, then \eqref{a}-\eqref{c} are derived from \eqref{v} and \eqref{b1}.
\ef
\begin{remark}
For given $t\in[0,\ve]$, it holds that
\begin{align}\label{x}
&\eta^{(0),-1}(x,t)=\beta^{(0)}_+(t)+(x-\phi(t))(\p_x\eta^{(0),-1})(\bar{x},t)>0\quad\text{for}~~x>\phi(t)+\f{6}{5}t^{\f{3}{2}},\no\\
&\eta^{(0),-1}(x,t)=\beta^{(0)}_-(t)+(x-\phi(t))(\p_x\eta^{(0),-1})(\bar{x},t)<0\quad\text{for}~~x<\phi(t)-\f{6}{5}t^{\f{3}{2}}.
\end{align}
\end{remark}

Note that by \eqref{eta0} and Lemma \ref{eta}, the solution of \eqref{w20} is
\beq\label{solu}
w^{(0)}_i(x,t)=w_i(\eta^{(0),-1}(x,t),0).
\eeq
From the properties of $\eta^{(0)}$ and $w_i(x,0)$, one has $w^{(0)}_i(x,t)\in C^2(\mathcal{B}_{\ve})$.
Differentiating both sides of \eqref{eta0} with respect to $x$,
we have that by \eqref{cz1} and \eqref{c},
\begin{align}
|\p_xw^{(0)}_i(x,t)|&=\Big|\f{w'_i(\eta^{(0),-1}(x,t),0)}{1+t\p_{w_i}\la_i(\overline{w}(\eta^{(0),-1}(x,t),0))w'_i(\eta^{(0),-1}(x,t),0)}\Big|\label{pxw}\\
&\leqslant\f{4}{3}(t^3+(x-\phi(t))^2)^{-\f{1}{3}}.\label{g1}
\end{align}
By taking the derivative with respect to $x$ twice on both sides of \eqref{eta0},
it is derived from \eqref{cz1}, \eqref{cz2} and \eqref{fm} that
\begin{align}
&|\p^2_xw^{(0)}_i(x,t)|\no\\
=&\Big|\f{w''_i(\eta^{(0),-1}(x,t),0)-t\p^2_{w_i}\la_i(\overline{w}(\eta^{(0),-1}(x,t),0))(w'_i(\eta^{(0),-1}(x,t),0))^3}{(1+t\p_{w_i}\la_i(\overline{w}(\eta^{(0),-1}(x,t),0))w'_i(\eta^{(0),-1}(x,t),0))^3}\Big|\no\\
\leqslant&10(t^3+(x-\phi(t))^2)^{-\f{5}{6}}.\label{g2}
\end{align}
Set
\begin{equation*}
[w^{(0)}_i(t)]=w^{(0)}_i(\phi(t)_-,t)-w^{(0)}_i(\phi(t)_+,t),\quad\langle w^{(0)}_i(t)\rangle=\f{1}{2}(w^{(0)}_i(\phi(t)_-,t)+w^{(0)}_i(\phi(t)_+,t)).
\end{equation*}
Together with \eqref{cz01}, \eqref{v} and \eqref{solu}, this yields
\begin{align}\label{w0j}
&|[w^{(0)}_i(t)]|\no\\
=&|w_i(\beta^{(0)}_-(t),0)-w_i(\beta^{(0)}_+(t),0)|\no\\
=&|-\beta^{(0)}_-(t)^{\f{1}{3}}+\beta^{(0)}_+(t)^{\f{1}{3}}+a_2\beta^{(0)}_-(t)^{\f{2}{3}}
-a_2\beta^{(0)}_+(t)^{\f{2}{3}}+O(\beta^{(0)}_-(t)-\beta^{(0)}_+(t))|\no\\
\leqslant&3t^{\f{1}{2}}.
\end{align}
Meanwhile, it follows from \eqref{jb}, \eqref{t} and \eqref{2-26} that
\begin{align}\label{wi0m}
&|\langle w^{(0)}_i\rangle|\no\\
\leqslant&\Big|\f{1}{2}(-\beta^{(0)}_-(t)^{\f{1}{3}}-\beta^{(0)}_+(t)^{\f{1}{3}})
+\f{1}{2}a_2(\beta^{(0)}_-(t)^{\f{2}{3}}+\beta^{(0)}_+(t)^{\f{2}{3}})+O(\beta^{(0)}_-(t)+\beta^{(0)}_+(t))\Big|\no\\
\leqslant& \f{2}{3}M^4t.
\end{align}
Therefore, in terms of \eqref{cz1}, \eqref{jb} and \eqref{g0},
we can arrive at
\begin{align}\label{wi0}
&\Big|\f{d}{dt}w_i(\beta^{(0)}_{\pm}(t),0)\pm\f{1}{2}t^{-\f{1}{2}}\Big|\no\\
=&\Big|w'_i(\beta^{(0)}_{\pm}(t),0)\f{d}{dt}\beta^{(0)}_{\pm}(t)\pm\f{1}{2}t^{-\f{1}{2}}\Big|\no\\
=&\Big|-\f{1}{3}(\beta^{(0)}_{\pm}(t))^{-\f{2}{3}}\f{\phi'(t)-\la_i(\overline{w}(\beta^{(0)}_{\pm}(t),0))}
{1+t\p_{w_i}\la_i(\overline{w}(\beta^{(0)}_{\pm}(t),0))w'_i(\beta^{(0)}_{\pm}(t),0)}
\pm\f{1}{2}t^{-\f{1}{2}}\Big|\no\\
&+\Big|\f{2}{3}a_2(\beta^{(0)}_{\pm}(t))^{-\f{1}{3}}\f{\phi'(t)-\la_i(\overline{w}(\beta^{(0)}_{\pm}(t),0))}
{1+t\p_{w_i}\la_i(\overline{w}(\beta^{(0)}_{\pm}(t),0))w'_i(\beta^{(0)}_{\pm}(t),0)}\Big|\no\\
\leqslant& M^6,
\end{align}
where the expression $\phi'(t)=\p_{\beta}\eta^{(0)}(\beta^{(0)}_{\pm}(t),t)\f{d}{dt}\beta^{(0)}_{\pm}(t)+\p_t\eta^{(0)}(\beta^{(0)}_{\pm}(t),t)$
is used.

\section {The construction of an approximate shock solution}\label{Sec-3}
Multiplying both sides of \eqref{1.3-1} by $l_l(w)$ for $l=1,2,\cdots,n$,
then \eqref{1.3-1} can be rewritten as
\begin{equation}\label{f}
\begin{cases}
\p_tw_j+\la_j(w)\p_xw_j+\ss_{k\neq i,j}p_{jk}(w)(\p_tw_k+\la_j\p_xw_k)=0,\quad j\neq i, \\
\p_tw_i+\la_i(w)\p_xw_i+\ss_{k\neq i}p_{ik}(w)(\p_tw_k+\la_i\p_xw_k)=0,
\end{cases}
\end{equation}
where $p_{lk}(w)=\f{l_{lk}(w)}{l_{ll}(w)}$ for $l=1,2,\cdots,n$ and $p_{lk}(0)=0$. Besides, the initial data of \eqref{f}
satisfies \eqref{2-19}.

Due to \eqref{utezheng}, for $|w_j|\lesssim\ve(j\neq i)$ and $|w_i|\lesssim\ve^{\f{1}{3}}$,
then there exists a large constant $M>0$ such that
\begin{equation}\label{3-1}
\la_{l+1}(w)-\la_l(w)\geqslant 4M^{-1},\quad l=1,\cdots,n-1.
\end{equation}
Set
\begin{align*}
&\Omega_+=\{(x,t)\in\mathcal{B}_{\ve};x>\phi(t)\};\quad\quad\quad\Omega_-=\{(x,t)\in\mathcal{B}_{\ve};x<\phi(t)\}.
\end{align*}
Then $\Omega_+\bigcup\Omega_-=\mathcal{B}_{\ve}$.

Define $w$  in $\Omega_{\pm}$ as $w_{\pm}(x,t)=(w_{1,\pm},\cdots, w_{n,\pm})(x,t)$, respectively.
In terms of the entropy condition \eqref{eq:3.28-1},
we decompose \eqref{f} with the initial data \eqref{2-19} in $\Omega_{\pm}$
into a coupled system of singular initial value problem and singular initial-boundary value problem
($w_i(x,0)$ has singularities). The initial data are derived from
the pre-shock and the boundary data are obtained from \eqref{RH}. Namely,
\beq\label{i}
\left\{
     \begin{array}{ll}
      \p_tw_{j,+}+\la_j(w_+)\p_xw_{j,+}+\ss_{k\neq i,j}p_{jk}(w_+)(\p_tw_{k,+}+\la_j(w_+)\p_xw_{k,+})=0,\quad j<i, \\
      \p_tw_{i,\pm}+\la_i(w_{\pm})\p_xw_{i,\pm}+\ss_{k\neq i}p_{ik}(w_{\pm})(\p_tw_{k,\pm}+\la_i(w_{\pm})\p_xw_{k,\pm})=0,\\
      \p_tw_{j,-}+\la_j(w_-)\p_xw_{j,-}+\ss_{k\neq i,j}p_{jk}(w_-)(\p_tw_{k,-}+\la_j(w_-)\p_xw_{k,-})=0,\quad j>i,\\
      w_{j,+}(x,t)|_{t=0}=0,\quad j<i,\\
      w_{i,\pm}(x,t)|_{t=0}=w_{i,\pm}(x,0),\\
      w_{j,-}(x,t)|_{t=0}=0,\quad j>i
     \end{array}
   \right.
\eeq
and
\beq\label{ib}
\left\{
     \begin{array}{ll}
     \p_tw_{j,-}+\la_j(w_-)\p_xw_{j,-}+\ss_{k\neq i,j}p_{jk}(w_-)(\p_tw_{k,-}+\la_j(w_-)\p_xw_{k,-})=0,\quad j<i, \\
     \p_tw_{j,+}+\la_j(w_+)\p_xw_{j,+}+\ss_{k\neq i,j}p_{jk}(w_+)(\p_tw_{k,+}+\la_j(w_+)\p_xw_{k,+})=0,\quad j>i, \\
      w_{j,-}(x,t)|_{t=0}=0,\quad w_{j,-}(x,t)|_{x=\phi(t)}=w_{j,-}(\phi(t)_-,t),\quad j<i,\\
      w_{j,+}(x,t)|_{t=0}=0,\quad w_{j,+}(x,t)|_{x=\phi(t)}=w_{j,+}(\phi(t)_+,t),\quad j>i.
     \end{array}
   \right.
\eeq
In the following discussion, we solve the systems \eqref{i} and \eqref{ib} together with \eqref{RH}.

\subsection{Iteration scheme for the approximate shock solution}
In this subsection, we discuss the iteration scheme for the approximate shock solution of \eqref{zhuyao}.
The iteration schemes are given as follows
\beq\label{im}
\left\{
     \begin{array}{ll}
    \p_tw^{(m+1)}_{j,+}+\la_j(w^{(m)}_+)\p_xw^{(m+1)}_{j,+}+\ss_{k\neq i,j}p_{jk}(w^{(m)}_+)(\p_tw^{(m)}_{k,+}+\la_j(w^{(m)}_+)\p_xw^{(m)}_{k,+})=0,\quad j<i, \\
      \p_tw^{(m+1)}_{i,\pm}+\la_i(w^{(m)}_{\pm})\p_xw^{(m+1)}_{i,\pm}+\ss_{k\neq i}p_{ik}(w^{(m)}_{\pm})(\p_tw^{(m)}_{k,\pm}+\la_i(w^{(m)}_{\pm})\p_xw^{(m)}_{k,\pm})=0,\\
      \p_tw^{(m+1)}_{j,-}+\la_j(w^{(m)}_-)\p_xw^{(m+1)}_{j,-}+\ss_{k\neq i,j}p_{jk}(w^{(m)}_-)(\p_tw^{(m)}_{k,-}+\la_j(w^{(m)}_-)\p_xw^{(m)}_{k,-})=0,\quad j>i,\\
      w^{(m+1)}_{j,+}(x,t)|_{t=0}=0,\quad j<i,\\
      w^{(m+1)}_{i,\pm}(x,t)|_{t=0}=w_{i,\pm}(x,0),\\
      w^{(m+1)}_{j,-}(x,t)|_{t=0}=0,\quad j>i
     \end{array}
   \right.
\eeq
and
\beq\label{ibm}
\left\{
     \begin{array}{ll}
     \p_tw^{(m+1)}_{j,-}+\la_j(w^{(m)}_-)\p_xw^{(m+1)}_{j,-}+\ss_{k\neq i,j}p_{jk}(w^{(m)}_-)(\p_tw^{(m)}_{k,-}+\la_j(w^{(m)}_-)\p_xw^{(m)}_{k,-})=0,\quad j<i,\\
     \p_tw^{(m+1)}_{j,+}+\la_j(w^{(m)}_+)\p_xw^{(m+1)}_{j,+}+\ss_{k\neq i,j}p_{jk}(w^{(m)}_+)(\p_tw^{(m)}_{k,+}+\la_j(w^{(m)}_+)\p_xw^{(m)}_{k,+})=0,\quad j>i, \\
     w^{(m+1)}_{j,-}(x,t)|_{t=0}=0,
     \quad w^{(m+1)}_{j,-}(x,t)|_{x=\phi(t)}=w^{(m+1)}_{j,-}(\phi(t)_-,t),\quad j<i,\\
     w^{(m+1)}_{j,+}(x,t)|_{t=0}=0,\quad w^{(m+1)}_{j,+}(x,t)|_{x=\phi(t)}=w^{(m+1)}_{j,+}(\phi(t)_+,t),\quad j>i.
     \end{array}
   \right.
\eeq

By the systems \eqref{im} and \eqref{ibm}, we will establish the stability of the iteration scheme
under iteration $m\longrightarrow m+1$.
Firstly, we define the iteration space as follows.

Let $M>1$ be a sufficiently large constant.
For all $(x,t)\in \mathcal{B}_{\ve}$, set
\begin{align*}
\chi^{(m)}_{\ve}=\{&(w^{(m)}_{1,\pm},\cdots,w^{(m)}_{i-1,\pm},w^{(m)}_{i,\pm},w^{(m)}_{i+1,\pm},\cdots,w^{(m)}_{n,\pm})\in C^1_{x,t}(\mathcal{B}_{\ve}):\no\\
&(w^{(m)}_{1,+},\cdots,w^{(m)}_{i-1,+},w^{(m)}_{i,\pm},w^{(m)}_{i+1,-},\cdots,w^{(m)}_{n,-})|_{t=0}=(0,\cdots,0,w_i(x,0),0,\cdots,0),\no\\
&\|(w^{(m)}_{1,\pm},\cdots,w^{(m)}_{i-1,\pm},w^{(m)}_{i,\pm},w^{(m)}_{i+1,\pm},\cdots,w^{(m)}_{n,\pm})\|_{\ve}\leqslant 1\},
\end{align*}
where the norm is defined as
\begin{align*}
&\|(w^{(m)}_{1,\pm},\cdots,w^{(m)}_{i-1,\pm},w^{(m)}_{i,\pm},w^{(m)}_{i+1,\pm},\cdots,w^{(m)}_{n,\pm})\|_{\ve}\no\\
=&\sup_{(x,t)\in\mathcal{B}_{\ve}}\{M^{-1}t^{-1}|w^{(m)}_{i,\pm}-w^{(0)}_{i,\pm}|, M^{-2}(t^3+(x-\phi(t))^2)^{\f{1}{6}}|\p_xw^{(m)}_{i,\pm}-\p_xw^{(0)}_{i,\pm}|,\no\\
&\quad\quad\quad~~M^{-3}(t^3+(x-\phi(t))^2)^{\f{1}{6}}|\p_tw^{(m)}_{i,\pm}-\p_tw^{(0)}_{i,\pm}|,
M^{-2}t^{-1}(t^3+(x-\phi(t))^2)^{-\f{1}{6}}|w^{(m)}_{j,\pm}|,\no\\
&\quad\quad\quad~~ M^{-4}(t^3+(x-\phi(t))^2)^{-\f{1}{6}}|\p_xw^{(m)}_{j,\pm}|,
M^{-5}(t^3+(x-\phi(t))^2)^{-\f{1}{6}}|\p_tw^{(m)}_{j,\pm}|\}.
\end{align*}
For convenience, we rewrite $\|(w^{(m)}_{1,\pm},\cdots,w^{(m)}_{i-1,\pm},w^{(m)}_{i,\pm},
w^{(m)}_{i+1,\pm},\cdots,w^{(m)}_{n,\pm})\|_{\ve}\leqslant 1$ as
\begin{align}
&|w^{(m)}_{i,\pm}(x,t)-w^{(0)}_{i,\pm}(x,t)|\leqslant Mt,\label{wm-w0}\\
&|\p_xw^{(m)}_{i,\pm}(x,t)-\p_xw^{(0)}_{i,\pm}(x,t)|\leqslant M^2(t^3+(x-\phi(t))^2)^{-\f{1}{6}},\label{pxwm}\\
&|\p_tw^{(m)}_{i,\pm}(x,t)-\p_tw^{(0)}_{i,\pm}(x,t)|\leqslant M^3(t^3+(x-\phi(t))^2)^{-\f{1}{6}},\label{ptwm}\\
&|w^{(m)}_{j,\pm}(x,t)|\leqslant M^2t(t^3+(x-\phi(t))^2)^{\f{1}{6}},\quad j\neq i,\label{wjm}\\
&|\p_xw^{(m)}_{j,\pm}(x,t)|\leqslant M^4(t^3+(x-\phi(t))^2)^{\f{1}{6}},\quad j\neq i,\label{pxwjm}\\
&|\p_tw^{(m)}_{j,\pm}(x,t)|\leqslant M^5(t^3+(x-\phi(t))^2)^{\f{1}{6}},\quad j\neq i.\label{ptwjm}
\end{align}

In addition, define the $i$-th characteristics $\eta^{(m)}(\beta,s)$ and
the $j$-th characteristics $\psi^{(m)}_{t,j}(x,s)$ $(j\neq i)$ in  \eqref{im} and \eqref{ibm} as
\begin{align}
  &\p_s\eta^{(m)}(\beta,s)=\la_i(w(\eta^{(m)}(\beta,s),s)),~~~~~~~~~~~~\eta^{(m)}(\beta,0)=\beta,\label{etam}\\
  &\p_s\psi^{(m)}_{t,j}(x,s)=\la_j(w(\psi^{(m)}_{t,j}(x,s),s)),~~~~~~~~~~\psi^{(m)}_{t,j}(x,t)=x,\quad j\neq i,\label{psi}
\end{align}
where $s\in[0,t]$. Note that $\eta^{(m)}$ is assigned an initial data at time $0$,
while $\psi^{(m)}_{t,j}$ is given a specified terminal data at time $t$.

\begin{lemma}\label{etam0}
Assume that $(w^{(m)}_{i,\pm}, w^{(m)}_{j,\pm})$ satisfies \eqref{wm-w0}-\eqref{ptwjm},
and the approximate shock curve $x=\phi(t)$ fulfills \eqref{jb}. Then for
$t\in[0,\ve]$, $\eta^{(m)}(\beta,t)$  in \eqref{etam} admits the following properties:
\begin{itemize}
  \item There exists a largest $\beta^{(m)}_+(t)>0$ and a smallest $\beta^{(m)}_-(t)<0$ satisfying
\begin{equation}\label{3-14}
\phi(t)=\eta^{(m)}(\beta^{(m)}_{\pm}(t),t),\quad\quad \f{4}{5}t^{\f{3}{2}}<|\beta^{(m)}_{\pm}(t)|<\f{6}{5}t^{\f{3}{2}}.
\end{equation}
 \item $\eta^{(m)}(\cdot,t):\Upsilon_m(t):=[-M\ve,M\ve]\backslash[\beta^{(m)}_-(t),\beta^{(m)}_+(t)]\longrightarrow[-M\ve,M\ve]\backslash\{\phi(t)\}$
is a bijection, and its inverse map is denoted by
$\eta^{(m),-1}(\cdot,t):[-M\ve,M\ve]\backslash\{\phi(t)\}\longrightarrow\Upsilon_m(t)$.
 \item  For $\beta\in\Upsilon_m(t)$, it holds that
\begin{align}
 &\f{1}{3}\leqslant \p_{\beta}\eta^{(m)}(\beta,t)\leqslant \f{33}{20},\label{pxeta}\\
 &|\eta^{(m)}(\beta,t)-\eta^{(0)}(\beta,t)|\leqslant 8Mt^2, \label{2phin-phi0}\\
 &|\p_{\beta}\eta^{(m)}(\beta,t)-\p_{\beta}\eta^{(0)}(\beta,t)|\leqslant M^3t^{\f{1}{2}}\label{pxetam}
\end{align}
and
\beq\label{w0etan}
\int^t_0|\p_x w^{(0)}_i(\eta^{(m)}(\beta,t'),t')|d t'\leqslant\f{24}{55}.
\eeq
\end{itemize}
\end{lemma}
\pf
It follows from \eqref{itz}, \eqref{g1} and \eqref{wm-w0}-\eqref{ptwjm} that
\begin{align}
&|\la^{(m)}_i(x,t)-\la^{(0)}_i(x,t)|\leqslant \f{3}{2}Mt,\label{lam}\\
&|\p_x\la^{(m)}_i(x,t)-\p_x\la^{(0)}_i(x,t)|\leqslant \f{11}{10}M^2(t^3+(x-\phi(t))^2)^{-\f{1}{6}}+M^7.\label{pxlam}
\end{align}
Introducing $q(\beta,t)$ with $q(\beta,0)=\beta$ such that $\eta^{(0)}(q(\beta,t),t)=\eta^{(m)}(\beta,t)$.
This yields
\begin{equation}\label{q}
q(\beta,t)=\eta^{(0),-1}(\eta^{(m)}(\beta,t),t)\in\Upsilon_0(t).
\end{equation}
In addition, from \eqref{a}, \eqref{lam}, \eqref{q} and the transport equation
$\p_t\eta^{(0),-1}+\la^{(0)}_i\p_x\eta^{(0),-1}=0$, we have
\begin{align}\label{dq}
 |\eta^{(0),-1}(\eta^{(m)}(\beta,t),t) & -\beta|=|q(\beta,t)-q(\beta,0)|=\Big|\int^t_0\p_{t'}q(\beta,t')dt'\Big| \no\\
 &=\Big|\int^t_0(\la^{(m)}_i-\la^{(0)}_i)\circ\eta^{(m)}(\p_x\eta^{(0),-1})\circ\eta^{(m)}dt'\Big|\leqslant 3Mt^2.
\end{align}
Due to $\phi(t)=\eta^{(0)}(\beta^{(0)}_{\pm}(t),t)=\eta^{(m)}(\beta^{(m)}_{\pm}(t),t)$,
together with \eqref{2-30} and \eqref{dq}, one then has
\begin{equation*}
  \f{4}{5}t^{\f{3}{2}}<|\beta^{(m)}_{\pm}(t)|<\f{6}{5}t^{\f{3}{2}}.
\end{equation*}
Note that we have from \eqref{etam},
\begin{equation}\label{3-23}
\p_{\beta}\eta^{(m)}(\beta,t)=e^{\int^t_0\p_x\la^{(m)}_i\circ\eta^{(m)}(\beta,t')d t'}.
\end{equation}
Since it follows from \eqref{cz1}, \eqref{czt}, \eqref{fm}, \eqref{pxw}, \eqref{q} and \eqref{dq} that
\begin{align*}
  &\int^t_0|\p_xw^{(0)}_i(\eta^{(m)}(\beta,t'),t')|dt' \no\\ =&\int^t_0\f{|w'_i(q(\beta,t'),0)|}{1+t'\p_{w_i}\la_i(\overline{w}(q(\beta,t'),0))w'_i(q(\beta,t'),0)}dt'\no\\
   \leqslant&\int^t_0\f{|w'_i(q(\beta,t'),0)|}{1+\f{11t'}{10}w'_i(q(\beta,t'),0)}dt'\no\\
   \leqslant&\f{10}{11}\int^t_0\f{\f{11}{10}|w'_i(\beta,0)|}{1+\f{11t'}{10}w'_i(\beta,0)}dt'
   +5M\int^t_0(t')^2\sup_{|\bar{\beta}-\beta|\leqslant3M(t')^2}\f{|w''_i(\bar{\beta},0)|}{(1+\f{11t'}{10}w'_i(\bar{\beta},0))^2}dt'\no\\
   \leqslant&\f{10}{11}\f{|w'_i(\beta,0)|}{w'_i(\beta,0)}\ln(1+\f{11t}{10}w'_i(\beta,0))+30Mt^{\f{1}{2}}\no\\
   \leqslant&\f{24}{55},
\end{align*}
where last inequality is derived by \eqref{cz1},
\begin{equation*}
\begin{cases}{}
  -\f{19}{50}\leqslant\f{11t}{10}w'_i(\beta,0)\leqslant0,~~~~~~~~~~~~~~~~~~~~ \f{19}{20}t^{\f{3}{2}}\leqslant|\beta|\leqslant t, \\[2mm]
   \f{11t}{10}|w'_i(\beta,0)|\leqslant\ve^{\f{1}{3}},~~~~~~~~~~~~~~~~~~~~~~~~~~~~~~|\beta|\geqslant t,
\end{cases}
\end{equation*}
and the fact of $\f{10}{11}\mathrm{sgn}(r)\ln(1+r)\leqslant\ln(\f{309}{200})\leqslant\f{109}{250}$ for all $r\in(-\f{19}{50},\ve^{\f{1}{3}})$,
then we obtain that by $\f{9}{10}\leqslant\p_{w_i}\la_i(\overline{w}(\beta,0))\leqslant\f{11}{10}$ and \eqref{w0etan},
\begin{equation}\label{3-24}
\int^t_0|\p_x\la^{(0)}_i\circ\eta^{(m)}|dt'\leqslant \f{12}{25}.
\end{equation}
In addition, \eqref{3-23} can be treated from
\eqref{pxlam} and \eqref{3-24} that
\begin{align}\label{3-26}
  \f{1}{3}\leqslant\p_{\beta}\eta^{(m)}(\beta,t) & =e^{\int^t_0\p_x\la^{(m)}_i\circ\eta^{(m)}(\beta,t')dt'} \no\\
   &= e^{\int^t_0\p_x\la^{(0)}_i\circ\eta^{(m)}(\beta,t')dt'}\cdot e^{\int^t_0(\p_x\la^{(m)}_i-\p_x\la^{(0)}_i)\circ\eta^{(m)}(\beta,t')dt'}\leqslant\f{33}{20},
   \end{align}
which means that the map $\eta^{(m)}(\cdot,t):\Upsilon_m(t)\longrightarrow[-M\ve,M\ve]\backslash\{\phi(t)\}$ is an injection.
Besides, the map $\eta^{(m)}$ is a surjection. Indeed, for $s\in[0,t)$, by \eqref{jb}, \eqref{v}, \eqref{g0}, \eqref{phixj}, \eqref{3-14} and \eqref{lam},
we arrive at
\begin{align}\label{phixn}
  \phi(s)-\eta^{(m)}(\beta^{(m)}_-(t),s)
  =&-\int^t_s(\phi'(t')-\la_i^{(m)}(w\circ\eta^{(m)}(\beta^{(m)}_-(t),t')))dt' \no\\
   =&\int^t_s[\la_i^{(0)}-\phi'(t')+(\la_i^{(m)}-\la_i^{(0)})]\circ\eta^{(m)}dt'\no\\
   \geqslant&\int^t_s\la_i(\overline{w}(q(\beta^{(m)}_-(t),t'),0))dt'-3M^4(t^2-s^2)\no\\
  \geqslant &-(\beta^{(m)}_-(t))^{\f{1}{3}}(t-s)-M^2(t^{\f{5}{3}}-s^{\f{5}{3}})-4M^4(t^2-s^2)\no\\
    \geqslant&\f{3}{5}t^{\f{1}{2}}(t-s).
\end{align}
This yields that by \eqref{b1}, \eqref{q} and \eqref{dq},
\begin{align*}
  |\eta^{(m)}(\beta,t)-\eta^{(0)}(\beta,t) |& =|\eta^{(0)}(q(\beta,t),t)-\eta^{(0)}(\beta,t)|\no\\
  &=|\p_{\beta}\eta^{(0)}(\overline{\beta})(q(\beta,t)-q(\beta,0))|\leqslant 6Mt^2,
\end{align*}
where $\overline{\beta}$ lies between $\beta$ and $q(\beta,t)$.

For $\beta\in\Upsilon_m(s)$ and  $s\in[0,t]$, collecting \eqref{jb}, \eqref{g0}, \eqref{2phin-phi0}, \eqref{lam} and \eqref{q},
one has
\begin{align}\label{s}
\sup_{s\in[0,t]}|\eta^{(m)}(\beta,s)-\phi(s)|
\leqslant&\sup_{s\in[0,t]}|\eta^{(m)}(\beta,s)-\eta^{(0)}(\beta,s)|+\sup_{s\in[0,t]}|\eta^{(0)}(\beta,s)-\phi(s)|\no\\
\leqslant&6Ms^2+\sup_{s\in[0,t]}|\beta+\int^s_0(\la^{(0)}_i\circ\eta^{(0)}-\phi'(t'))dt'|\no\\
\leqslant&6Ms^2+|\beta|+\sup_{s\in[0,t]}\int^s_0(|\la^{(0)}_i|+|\phi'(t')|)\circ\eta^{(0)}dt'\no\\
\leqslant&|\beta|+\sup_{s\in[0,t]}2s|\beta|^{\f{1}{3}}+2M^4s^2\no\\
\leqslant&2M\ve.
\end{align}
On the other hand, it is derived from \eqref{eta0} and \eqref{etam} that
\begin{align}\label{pbetaphin0}
  &\p_t(\p_{\beta}\eta^{(m)}(\beta,t)-\p_{\beta}\eta^{(0)}(\beta,t))\no\\
=&\p_x\la^{(0)}_i\circ\eta^{(m)}(\beta,t)(\p_{\beta}\eta^{(m)}(\beta,t)-\p_{\beta}\eta^{(0)}(\beta,t))\no\\
&+(\p_x\la^{(m)}_i\circ\eta^{(m)}(\beta,t)-\p_x\la^{(0)}_i\circ\eta^{(m)}(\beta,t))\p_{\beta}\eta^{(m)}(\beta,t)\no\\
   &+(\p_x\la^{(0)}_i\circ\eta^{(m)}(\beta,t)-\p_x\la^{(0)}_i\circ\eta^{(0)}(\beta,t))\p_{\beta}\eta^{(0)}(\beta,t).
\end{align}
Based on \eqref{czt}, \eqref{fm}, \eqref{b1}, \eqref{pxw}, \eqref{pxeta}, \eqref{pxlam} and \eqref{dq}, we have
\begin{equation}\label{3-28}
|(\p_x\la^{(m)}_i\circ\eta^{(m)}(\beta,t)-\p_x\la^{(0)}_i\circ\eta^{(m)}(\beta,t))\p_{\beta}\eta^{(m)}(\beta,t)|\leqslant 4 M^2 t^{-\f{1}{2}}
\end{equation}
and
\begin{align}\label{3-29}
 &|(\p_x\la^{(0)}_i\circ\eta^{(m)}(\beta,t)-\p_x\la^{(0)}_i\circ\eta^{(0)}(\beta,t))\p_{\beta}\eta^{(0)}(\beta,t) |\no\\
 \leqslant& \f{11}{10}\Big|\p_xw^{(0)}_i\circ\eta^{(m)}(\beta,t)-\p_xw^{(0)}_i\circ\eta^{(0)}(\beta,t)\Big||\p_{\beta}\eta^{(0)}(\beta,t)| \no\\
   \leqslant& 2\Big|\f{w'_i(q(\beta,t),0)}{1+t\p_{w_i}\la_i(\overline{w}(q(\beta,t),0))w'_i(q(\beta,t),0)}
   -\f{w'_i(\beta,0)}{1+t\p_{w_i}\la_i(\overline{w}(\beta,0))w'_i(\beta,0)}\Big||\p_{\beta}\eta^{(0)}(\beta,t)|\no\\
   \leqslant&2\Bigg|\f{w''_i(\overline{\beta},0)}{(1+\f{11t}{10}w'_i(\overline{\beta},0))^2}\Bigg||q(\beta,t)
   -q(\beta,0)||\p_{\beta}\eta^{(0)}(\beta,t)|\leqslant 25Mt^{-\f{1}{2}}.
\end{align}
Substituting \eqref{3-24}, \eqref{3-28} and \eqref{3-29}
into \eqref{pbetaphin0} and applying Growall's inequality, \eqref{pxetam} is proved.
\ef

Along the $j$-th characteristics  $\psi^{(m)}_{t,j}(x,s)$ ($j\not=i$), we treat the integral of $\p_xw^{(0)}_i$.

\begin{lemma}\label{haosan}
For $t'\in[s,t]$ and $(x,t)\in{\mathcal{B}_{\ve}}$, if there exists a point $\gamma(x,t')=\eta^{(0),-1}(\psi^{(m)}_{t,j}(x,t'),t')\in\Upsilon_0(t')$ such that $\psi^{(m)}_{t,j}(x,t')=\eta^{(0)}(\gamma(x,t'),t')$, then it holds that

\begin{itemize}
 \item  when $\gamma(x,t')\geqslant\f{6}{5}(t')^{\f{3}{2}}$ for $j>i$ or $\gamma(x,t')\leqslant-\f{6}{5}(t')^{\f{3}{2}}$ for $j<i$,
 one has
 \beq\label{pw0psin}
\int^t_{s}|\p_x w^{(0)}_i(\psi^{(m)}_{t,j}(x,t'),t')|d t'\leqslant\f{40}{|\la_j-\la_i|(0)}(t^3+(x-\phi(t))^2)^{\f{1}{6}};
\eeq
 \item when $\gamma(x,t')\geqslant\f{6}{5}(t')^{\f{3}{2}}$ for $j<i$ or $\gamma(x,t')\leqslant-\f{6}{5}(t')^{\f{3}{2}}$ for $j>i$,
 one has
\begin{equation*}
\int^t_{s}|\p_x w^{(0)}_i(\psi^{(m)}_{t,j}(x,t'),t')|d t'\leqslant\f{40}{|\la_j-\la_i|(0)}(s^3+(\psi^{(m)}_{t,j}(x,s)-\phi(s))^2)^{\f{1}{6}}.
\end{equation*}
\end{itemize}
\end{lemma}
\pf
It follows from \eqref{psi} and $\p_{t'}\eta^{(0),-1}+\la^{(0)}_i\p_x\eta^{(0),-1}=0$ that
\begin{align}\label{beta1d}
\p_{t'}\gamma(x,t')&=(\la^{(m)}_j(w(\psi^{(m)}_{t,j}(x,t'),t'))-\la^{(0)}_i(\overline{w}(\psi^{(m)}_{t,j}(x,t'),t')))
\p_x\eta^{(0),-1}(\psi^{(m)}_{t,j}(x,t'),t').
\end{align}
This, together with  \eqref{cz01}, \eqref{eta0}, \eqref{a}, \eqref{wm-w0}, \eqref{wjm} and \eqref{beta1d}, yields
\begin{align}
 &\p_{t'}\gamma(x,t')\geqslant\f{5}{16}(\la_j-\la_i)(0)>0,\quad~~~ j>i,\label{j}\\
 &\p_{t'}\gamma(x,t')\leqslant-\f{5}{16}(\la_i-\la_j)(0)<0,\quad j<i.\label{k}
\end{align}
Therefore, along $\psi^{(m)}_{t,j}$ for $s\in[0,t)$, by \eqref{cz1}, \eqref{fm}, \eqref{c}, \eqref{pxw}, \eqref{3-1}, \eqref{j} and \eqref{k}, we have
\begin{align}\label{l}
  &\int^t_{s}|\p_xw^{(0)}_i(\psi^{(m)}_{t,j}(x,t'),t')|dt'\no\\
  =& \int^t_{s}\f{|w'_i(\gamma(x,t'),0)|}{1+t'\p_{w_i}\la_i(\overline{w}(\gamma(x,t'),0))w'_i(\gamma(x,t'),0)}dt'\no\\
  \leqslant& 3\int^t_{s}(\gamma(x,t'))^{-\f{2}{3}}dt'\no\\
   \leqslant&\f{48}{5(\la_j-\la_i)(0)}\int^t_{s}\p_{t'}\gamma(x,t')(\gamma(x,t'))^{-\f{2}{3}}dt'\no\\
   =&
   \begin{cases}
   \f{144}{5(\la_j-\la_i)(0)}(\gamma(x,t)^{\f{1}{3}}-\gamma(x,s)^{\f{1}{3}}),\quad j>i,\\[2mm]
   \f{144}{5(\la_i-\la_j)(0)}(\gamma(x,s)^{\f{1}{3}}-\gamma(x,t)^{\f{1}{3}}),\quad j<i,
   \end{cases}
   \no\\
   \leqslant&
   \begin{cases}
        \f{32}{|\la_j-\la_i|(0)}|\gamma(x,t)|^{\f{1}{3}},  ~~\gamma(x,t')\geqslant\f{6}{5}(t')^{\f{3}{2}},j>i; \gamma(x,t')\leqslant-\f{6}{5}(t')^{\f{3}{2}}, j<i,\\[2mm]
        \f{32}{|\la_j-\la_i|(0)}|\gamma(x,s)|^{\f{1}{3}},
        ~~\gamma(x,t')\geqslant\f{6}{5}(t')^{\f{3}{2}}, j<i; \gamma(x,t')\leqslant-\f{6}{5}(t')^{\f{3}{2}}, j>i,
   \end{cases}
\end{align}
where the last inequality is due to \eqref{c}, \eqref{x}, \eqref{j} and \eqref{k}. Indeed,
note that $\gamma(x,t')^{\f{1}{3}}$ is a strictly increasing function with respect to $t'$ for $j>i$, which is derived from \eqref{j}.
Therefore, when $\gamma(x,t')\geqslant\f{6}{5}(t')^{\f{3}{2}}>0$ for all $t'\in[s,t]$,
$|\gamma(x,t)|^{\f{1}{3}}$ is the maximum value; when $\gamma(x,t')\leqslant-\f{6}{5}(t')^{\f{3}{2}}<0$ for all $t'\in[s,t]$,
the maximum value of  $|\gamma(x,t')|$ is $|\gamma(x,s)|^{\f{1}{3}}$. Analogous properties hold
for the situation of $j<i$.
\ef

Based on Lemma \ref{haosan}, we further have
\begin{lemma}\label{bt4.3}
For $\al=-\f{1}{3},-\f{1}{6},\f{1}{6},\f{1}{3},\f{2}{3}$, the following estimates hold under the conditions of Lemma \ref{haosan},
\begin{equation}\label{h}
\int^t_0|((t')^3+(\psi^{(m)}_{t,j}(x,t')-\phi(t'))^2)^{\al}|dt'\leqslant
\f{5M}{(2\al+1)|\la_j-\la_i|(0)}(t^3+(x-\phi(t))^2)^{\f{2\al+1}{3}}.
\end{equation}
\end{lemma}
\pf
From \eqref{c} and Lemma \ref{haosan}, $|((t')^3+(\psi^{(m)}_{t,j}(x,t')-\phi(t'))^2)^{\al}|$ is equivalent to $|\gamma(x,t')|^{2\al}$.
As in \eqref{l}, we can arrive at
\begin{equation*}
\begin{aligned}
&\int^t_0|((t')^3+(\psi^{(m)}_{t,j}(x,t')-\phi(t'))^2)^{\al}|dt'\\
\leqslant&\begin{cases}
     \f{5}{(2\al+1)|\la_j-\la_i|(0)}(t^3+(x-\phi(t))^2)^{\f{2\al+1}{2}},~~\gamma(x,t')\geqslant\f{6}{5}(t')^{\f{3}{2}}, j>i; \gamma(x,t')\leqslant-\f{6}{5}(t')^{\f{3}{2}}, j<i,\\[2mm]
     \f{5}{(2\al+1)|\la_j-\la_i|(0)}|\psi^{(m)}_{t,j}(x,0)|^{2\al+1},~~~~~~~~~~~~\gamma(x,t')\geqslant\f{6}{5}(t')^{\f{3}{2}}, j<i; \gamma(x,t')\leqslant-\f{6}{5}(t')^{\f{3}{2}}, j>i.
    \end{cases}
\end{aligned}
\end{equation*}
In addition, collecting \eqref{cz01}, \eqref{jb}, \eqref{3-1}, \eqref{wm-w0}, \eqref{wjm} and \eqref{psi} yields
\begin{align*}
    |\psi^{(m)}_{t,j}(x,0)|=&|\psi^{(m)}_{t,j}(x,t)-\phi(t)-\int^t_0(\la^{(m)}_j\circ\psi^{(m)}_{t,j}-\phi'(t'))dt'|\no\\
    \leqslant&|x-\phi(t)|+2|\la_j-\la_i|(0)t.
    \end{align*}
Thus,
\begin{align*}
        |\psi^{(m)}_{t,j}(x,0)|^{2\al+1}
        &\leqslant
        \begin{cases}
             M(t^3+(x-\phi(t))^2)^{\f{2\al+1}{2}},~~|x-\phi(t)|\geqslant 2|\la_j-\la_i|(0)t,  \\[2mm]
Mt^{2\al+1},~~\quad\quad\quad\quad\quad\quad\quad~|x-\phi(t)|\leqslant2|\la_j-\la_i|(0)t.
       \end{cases}\no\\
&\leqslant M(t^3+(x-\phi(t))^2)^{\f{2\al+1}{3}}.
    \end{align*}
Then \eqref{h} is obtained.
\ef

We turn our attention to the existence and uniqueness of $\psi^{(m)}_{t,j}(j\neq i)$ in \eqref{psi}.
Due to $|\la^{(m)}_j|\leqslant M$, \eqref{g1}, \eqref{pxwm} and \eqref{pxwjm}, one has
\begin{align}\label{3-39}
|\p_x\la^{(m)}_j&(w(\psi^{(m)}_{t,j}(x,t'),t'))|\leqslant\f{8}{3}M((t')^3+(\psi^{(m)}_{t,j}(x,t')-\phi(t'))^2)^{-\f{1}{3}}\no\\
&~~~~+M^3((t')^3+(\psi^{(m)}_{t,j}(x,t')-\phi(t'))^2)^{-\f{1}{6}}
+M^5((t')^3+(\psi^{(m)}_{t,j}(x,t')-\phi(t'))^2)^{\f{1}{6}}.
\end{align}
Then there exists a unique solution to \eqref{psi} for $(x,t)\in\mathcal{B}_{\ve}$.
Taking the first order derivative with respect to $x$ on both sides of \eqref{psi} and
integrating over $[s,t](s\in[0,t])$, we arrive at
\beq\label{pxpsi}
\p_x\psi^{(m)}_{t,j}(x,s)=e^{-\int^t_s\p_x\la^{(m)}_j(\psi^{(m)}_{t,j}(x,t'),t')dt'}.
\eeq
 In addition, by making use of \eqref{h} in Lemma \ref{bt4.3} and \eqref{3-39}, one can obtain a more precise estimate of \eqref{pxpsi}.
 In fact, for $|x-\phi(t)|\leqslant\ve^{\f{1}{2}}$, it holds that
\beq\label{suppsi}
\sup_{s\in[0,t]}|\p_x\psi^{(m)}_{t,j}(x,s)-1|\leqslant M\ve^{\f{1}{9}}.
\eeq

\subsection{Stability of the iteration scheme}
In this subsection, we investigate the stability of the iteration scheme.
\begin{lemma}\label{diedai}
For all $m\geqslant0$ and sufficiently large $M>0$, the map $\chi^{(m)}_{\ve}\longrightarrow\chi^{(m+1)}_{\ve}$ corresponding to
\begin{equation*}
(w^{(m)}_{1,\pm},\cdots,w^{(m)}_{i,\pm},\cdots,w^{(m)}_{n,\pm})\longmapsto(w^{(m+1)}_{1,\pm},\cdots,w^{(m+1)}_{i,\pm},\cdots,w^{(m+1)}_{n,\pm})
\end{equation*}
is stable. Namely, $(w^{(m+1)}_{1,\pm},\cdots,w^{(m+1)}_{i,\pm},\cdots,w^{(m+1)}_{n,\pm})$ satisfies \eqref{wm-w0}-\eqref{ptwjm}.
\end{lemma}
\pf
In the procedure of proof, we will appply the 0-th iteration \eqref{w20} and the iteration systems \eqref{im}-\eqref{ibm}.

\noindent{\bf Part 1. Estimate of $w^{(m+1)}_{i,\pm}(x,t)$}

Note that
\beq\label{wm+1-w0}
     \begin{cases}
       \p_t(w^{(m+1)}_{i,+}-w^{(0)}_{i,+})+\la_i(w^{(m)}_+)\p_x(w^{(m+1)}_{i,+}-w^{(0)}_{i,+})
       +(\la_i(w^{(m)}_+)-\la_i(w^{(0)}_+))\p_xw^{(0)}_{i,+}  \\[2mm]
       +\ss_{k\neq i}p_{ik}(w^{(m)}_+)(\p_tw^{(m)}_{k,+}+\la_i(w^{(m)}_{+})\p_xw^{(m)}_{k,+})=0,\\[4mm]
       w^{(m+1)}_{i,+}(x,0)-w^{(0)}_{i,+}(x,0)=0.
     \end{cases}
\eeq
Along the characteristics $x=\eta^{(m)}(\beta,t)$ for $|\beta|\leqslant M\ve<M^{-100}$, we obtain
\begin{align*}
(w^{(m+1)}_{i,+}-w^{(0)}_{i,+})\circ\eta^{(m)}(\beta,t)
=&-\int^t_{0}(\la_i(w^{(m)}_+)-\la_i(w^{(0)}_+))\circ\eta^{(m)}(\beta,t')\p_xw^{(0)}_{i,+}\circ\eta^{(m)}(\beta,t')dt' \no\\
& -\int^t_0\ss_{k\neq i}p_{ik}(w^{(m)}_+)\circ\eta^{(m)}(\beta,t')\Big(\f{d}{dt'}w^{(m)}_{k,+}\circ\eta^{(m)}(\beta,t')\Big)dt'.
\end{align*}
It is derived from \eqref{pxwjm}, \eqref{ptwjm}, \eqref{w0etan}, \eqref{lam} and \eqref{s} that
\begin{align*}
&|w^{(m+1)}_{i,+}(x,t)-w^{(0)}_{i,+}(x,t)|\no\\
\leqslant&\f{36}{55}Mt+M^{12}|2M\ve|^{\f{1}{3}}\int^t_0((t')^3+(\eta^{(m)}(\beta,t')-\phi(t'))^2)^{\f{1}{6}}dt'\no\\
\leqslant&\f{36}{55}Mt+2M^{13}\ve^{\f{1}{3}}(|2M\ve|^{\f{1}{3}}t+t^{\f{3}{2}})\no\\
\leqslant&\f{3}{4}Mt.
\end{align*}
where we have used the fact of $|p_{ik}(w^{(m)}_+)\circ\eta^{(m)}(\beta,t)|\leqslant M^6|2M\ve|^{\f{1}{3}}$.
In fact, by \eqref{cz01}, \eqref{wm-w0}, \eqref{wjm} and \eqref{s},
it holds that
\begin{align}\label{3-45}
&|p_{ik}(w^{(m)}_+)\circ\eta^{(m)}(\beta,t)|
=|p_{ik}(w^{(m)}_+)\circ\eta^{(m)}(\beta,t)-p_{ik}(0)|\no\\
\leqslant&|\ss_{k\neq i}\p_{w_{k,+}}p_{ik}(\widetilde{w}^{(m)}_+)w^{(m)}_{k,+}\circ\eta^{(m)}(\beta,t)+\p_{w_{i,+}}p_{ik}(\widetilde{w}^{(m)}_+)w^{(m)}_{i,+}\circ\eta^{(m)}(\beta,t)|\no\\
\leqslant& M^6|2M\ve|^{\f{1}{3}},
\end{align}
where $\widetilde{w}^{(m)}_+$ lies between $0$ and $w^{(m)}_+$.

Taking the first order derivative with respect to $x$ on both sides of \eqref{wm+1-w0}, one has
\beq
     \begin{cases}\label{pxwm+1-pxw0}
       \p_t(\p_xw^{(m+1)}_{i,+}-\p_xw^{(0)}_{i,+})+\la_i(w^{(m)}_+)\p^2_x(w^{(m+1)}_{i,+}-w^{(0)}_{i,+})
       +\p_x\la_i(w^{(m)}_+)(\p_xw^{(m+1)}_{i,+}-\p_xw^{(0)}_{i,+})\\[2mm]
       +(\p_x\la_i(w^{(m)}_+)-\p_x\la_i(w^{(0)}_+))\p_xw^{(0)}_{i,+}+(\la_i(w^{(m)}_+)-\la_i(w^{(0)}_+))\p^2_xw^{(0)}_{i,+}  \\[2mm]
       +\ss_{k\neq i}p_{ik}(w^{(m)}_+)(\p_t\p_xw^{(m)}_{k,+}+\p_x\la_i(w^{(m)}_{+})\p_xw^{(m)}_{k,+}
       +\la_i(w^{(m)}_{+})\p^2_xw^{(m)}_{k,+})\\[2mm]
       +\ss_{k\neq i}\p_xp_{ik}(w^{(m)}_+)(\p_tw^{(m)}_{k,+}+\la_i(w^{(m)}_{+})\p_xw^{(m)}_{k,+})=0,\\[2mm]
       \p_xw^{(m+1)}_{i,+}(x,0)-\p_xw^{(0)}_{i,+}(x,0)=0.
     \end{cases}
\eeq
Along the characteristics $x=\eta^{(m)}(\beta,t)$, we can rewrite the equation in \eqref{pxwm+1-pxw0} as
\begin{align}\label{ddtw}
&\f{d}{dt}((\p_xw^{(m+1)}_{i,+}-\p_xw^{(0)}_{i,+})\circ\eta^{(m)})+\p_x\la_i(w^{(m)}_+)\circ\eta^{(m)}(\p_xw^{(m+1)}_{i,+}
-\p_xw^{(0)}_{i,+})\circ\eta^{(m)}\no\\
=&-(\p_x\la_i(w^{(m)}_+)-\p_x\la_i(w^{(0)}_+))\circ\eta^{(m)}\p_xw^{(0)}_{i,+}\circ\eta^{(m)}
-(\la_i(w^{(m)}_+)-\la_i(w^{(0)}_+))\circ\eta^{(m)}\p^2_xw^{(0)}_{i,+}\circ\eta^{(m)}\no\\
&-\ss_{k\neq i}\f{d}{dt}(p_{ik}(w^{(m)}_+)\circ\eta^{(m)}\p_xw^{(m)}_{k,+}\circ\eta^{(m)}\p_{\beta}\eta^{(m)})(\p_{\beta}\eta^{(m)})^{-1}\no\\
&+\ss_{k\neq i}(\p_tp_{ik}(w^{(m)}_+)\p_xw^{(m)}_{k,+}-\p_tw^{(m)}_{k,+}\p_xp_{ik}(w^{(m)}_+))\circ\eta^{(m)}.
\end{align}
Then solving \eqref{ddtw} yields
\begin{align*}
&\p_xw^{(m+1)}_{i,+}(x,t)-\p_xw^{(0)}_{i,+}(x,t)\no\\
=&-\int^t_0e^{-\int^t_{t'}\p_x\la_i(w^{(m)}_+\circ\eta^{(m)})dt''}(\p_x\la_i(w^{(m)}_+)-\p_x\la_i(w^{(0)}_+))\circ\eta^{(m)}\p_xw^{(0)}_{i,+}\circ\eta^{(m)}dt'\no\\
&-\int^t_0e^{-\int^t_{t'}\p_x\la_i(w^{(m)}_+\circ\eta^{(m)})dt''}(\la_i(w^{(m)}_+)-\la_i(w^{(0)}_+))\circ\eta^{(m)}\p^2_xw^{(0)}_{i,+}\circ\eta^{(m)}dt'\no\\
&-\int^t_0e^{-\int^t_{t'}\p_x\la_i(w^{(m)}_+\circ\eta^{(m)})dt''}\ss_{k\neq i}\f{d}{dt'}(p_{ik}(w^{(m)}_+)\circ\eta^{(m)}\p_xw^{(m)}_{k,+}\circ\eta^{(m)}\p_{\beta}\eta^{(m)})(\p_{\beta}\eta^{(m)})^{-1}dt'\no\\
&+\int^t_0e^{-\int^t_{t'}\p_x\la_i(w^{(m)}_+\circ\eta^{(m)})dt''}\ss_{k\neq i}(\p_tp_{ik}(w^{(m)}_+)\p_xw^{(m)}_{k,+}-\p_tw^{(m)}_{k,+}\p_xp_{ik}(w^{(m)}_+))\circ\eta^{(m)}dt'\no\\
:=&I_1+I_2+I_3+I_4.
\end{align*}
Similarly to the treatment of \eqref{3-26}, we have
\beq\label{ef}
|e^{-\int^t_{t'}\p_x\la_i(w^{(m)}_+\circ\eta^{(m)})dt''}|\leqslant\f{33}{20}.
\eeq
It follows from \eqref{g2}, \eqref{lam}, \eqref{pxlam} and \eqref{ef} that
\begin{align}\label{3-50}
|I_1+I_2|\leqslant&\f{363}{200}M^2\int^t_0|\p_xw^{(0)}_{i,+}\circ\eta^{(m)}|((t')^3+(\eta^{(m)}(\beta,t')-\phi(t'))^2)^{-\f{1}{6}}dt'\no\\
&+30M\int^t_0t'((t')^3+(\eta^{(m)}(\beta,t')-\phi(t'))^2)^{-\f{5}{6}}dt'+\f{18}{25}M^7.
\end{align}
Denote by $x=\eta^{(m)}(\beta,t)$. As treated in \eqref{phixn}, for $|x|\geqslant\f{3}{4}t^{\f{3}{2}}$ and $t'\in[0,t]$, it is derived from \eqref{jb}, \eqref{v}, \eqref{g0}
and \eqref{lam} that
\begin{align}\label{3-51}
&\eta^{(m)}(\beta,t')-\phi(t')
=\eta^{(m)}(\beta,t)-\phi(t)-\int^t_{t'}(\la_i(w^{(m)}_+\circ\eta^{(m)})-\phi'(t''))dt''\no\\
&\geqslant x-\phi(t)-(t-t')(\la_i(\overline{w}(x,0)))-M^4(t^2-(t')^2)\no\\
&\geqslant x-\phi(t)+\f{\sqrt{3}}{2}t^{\f{1}{2}}(t-t')
\end{align}
and
\begin{align}\label{tt'}
(t')^3+(\eta^{(m)}(\beta,t')-\phi(t'))^2 &\geqslant(t')^3+(x-\phi(t))^2+\f{3}{4}t(t-t')^2\no\\
&\geqslant \f{5}{16}(t^3+(x-\phi(t))^2).
\end{align}
Thus, \eqref{3-50} admits the following estimate
\begin{align*}
|I_1+I_2|&\leqslant\Big(\f{363}{200}\f{24}{55}\Big(\f{5}{16}\Big)^{-\f{1}{6}}M^2+50M\Big)(t^3+(x-\phi(t))^2)^{-\f{1}{6}}+\f{4}{5}M^7\no\\
&\leqslant\Big(\f{97}{100}M^2+60M\Big)(t^3+(x-\phi(t))^2)^{-\f{1}{6}}.
\end{align*}
On the other hand,
\begin{align*}
I_3=&-\ss_{k\neq i}e^{-\int^t_{t'}\p_x\la_i(w^{(m)}_+\circ\eta^{(m)})dt''}p_{ik}(w^{(m)}_+)\circ\eta^{(m)}(\beta,t)\p_xw^{(m)}_{k,+}\circ\eta^{(m)}(\beta,t)\no\\
&+\ss_{k\neq i}\int^t_0p_{ik}(w^{(m)}_+)\circ\eta^{(m)}\p_xw^{(m)}_{k,+}\circ\eta^{(m)}\p_{\beta}\eta^{(m)}\f{d}{dt'}(e^{-\int^t_{t'}\p_x\la_i(w^{(m)}_+\circ\eta^{(m)})dt''}(\p_{\beta}\eta^{(m)})^{-1})dt',
\end{align*}
where
\begin{align*}
&\f{d}{dt'}(e^{-\int^t_{t'}\p_x\la_i(w^{(m)}_+\circ\eta^{(m)})dt''}(\p_{\beta}\eta^{(m)})^{-1})\no\\
=&\p_x\la_i(w^{(m)}_+\circ\eta^{(m)})(e^{-\int^t_{t'}\p_x\la_i(w^{(m)}_+\circ\eta^{(m)})dt''}(\p_{\beta}\eta^{(m)})^{-1})\no\\
&-e^{-\int^t_{t'}\p_x\la_i(w^{(m)}_+\circ\eta^{(m)})dt''}(\p_{\beta}\eta^{(m)})^{-1}\p_x\la_i(w^{(m)}_+\circ\eta^{(m)})\no\\
=&0.
\end{align*}
$I_4$ can be analogously treated.
Then it follows from \eqref{wm-w0}, \eqref{pxwm}, \eqref{ptwm}, \eqref{pxwjm}, \eqref{ptwjm}, \eqref{pxeta} and \eqref{s} that
\begin{equation*}
|I_3+I_4|\leqslant M.
\end{equation*}
Therefore, we arrive at
\begin{equation}\label{Sun-S1}
|\p_x w^{(m+1)}_{i,+}(x,t)-\p_x w^{(0)}_{i,+}(x,t)|\leqslant\f{99}{100}M^2(t^3+(x-\phi(t))^2)^{-\f{1}{6}}.
\end{equation}
Besides, together with \eqref{g1}, \eqref{wm-w0}, \eqref{pxwm}, \eqref{pxwjm}, \eqref{lam} and \eqref{wm+1-w0}, one has
\begin{equation}\label{Sun-S2}
|\p_t w^{(m+1)}_{i,+}(x,t)-\p_t w^{(0)}_{i,+}(x,t)|\leqslant\f{1}{2}M^3(t^3+(x-\phi(t))^2)^{-\f{1}{6}}.
\end{equation}
In the same way, we also obtain the analogous estimates for $w^{(m+1)}_{i,-}-w^{(0)}_{i,-}$
as in \eqref{Sun-S1} and \eqref{Sun-S2}.

\noindent{\bf Part 2. Estimate for $w^{(m+1)}_{j,\pm}(x,t)$}
\vskip 0.2 true cm

Note that for $j<i$,
we consider the following Cauchy problem
\beq\label{wjm+1C}
     \begin{cases}
     \p_tw^{(m+1)}_{j,+}+\la_j(w^{(m)}_+)\p_xw^{(m+1)}_{j,+}+\ss_{k\neq i,j}p_{jk}(w^{(m)}_+)(\p_tw^{(m)}_{k,+}+\la_j(w^{(m)}_+)\p_xw^{(m)}_{k,+})=0,\quad j<i, \\[2mm]
      w^{(m+1)}_{j,+}(x,t)|_{t=0}=0,\quad j<i.
     \end{cases}
\eeq
Along the characteristics $\psi^{(m)}_{t,j}(x,s)$ with $\psi^{(m)}_{t,j}(x,t)=x$ and $s\in[0,t]$, the equation in \eqref{wjm+1C} can be rewritten as
\begin{align}\label{wjm+1t}
&\f{d}{ds}(w^{(m+1)}_{j,+}(\psi^{(m)}_{t,j}(x,s),s)+\ss_{k\neq i,j}p_{jk}(w^{(m)}_+(\psi^{(m)}_{t,j}(x,s),s))w^{(m)}_{k,+}(\psi^{(m)}_{t,j}(x,s),s))\no\\
=&\ss_{k\neq i,j}\ss_{l}\p_{w_{l,+}}p_{jk}(w^{(m)}_+(\psi^{(m)}_{t,j}(x,s),s))\f{d}{ds}w^{(m)}_{l,+}(\psi^{(m)}_{t,j}(x,s),s)w^{(m)}_{k,+}(\psi^{(m)}_{t,j}(x,s),s).
\end{align}
Then integrating \eqref{wjm+1t} over $[0,t]$, together with the initial data in \eqref{wjm+1C}, yields
\begin{align}\label{wjm+1e}
&w^{(m+1)}_{j,+}(x,t)\no\\
=&-\ss_{k\neq i,j}p_{jk}(w^{(m)}_+(x,t))w^{(m)}_{k,+}(x,t)\no\\
&+\ss_{k\neq i,j}\int^t_0(\p_{w_{i,+}}p_{jk}(w^{(m)}_+)(\p_sw^{(m)}_{i,+}+\la_j(w^{(m)}_+)\p_xw^{(m)}_{i,+})w^{(m)}_{k,+})(\psi^{(m)}_{t,j}(x,s),s)ds\no\\
&+\ss_{k\neq i,j}\ss_{l\neq i}\int^t_0(\p_{w_{l,+}}p_{jk}(w^{(m)}_+)(\p_sw^{(m)}_{l,+}+\la_j(w^{(m)}_+)\p_xw^{(m)}_{l,+})w^{(m)}_{k,+})(\psi^{(m)}_{t,j}(x,s),s)ds.
\end{align}
Note that by \eqref{g1}, \eqref{wm-w0}-\eqref{ptwjm}, \eqref{psi}, \eqref{pw0psin} and \eqref{h} for $\al=-\f{1}{6}$,
we have
\begin{align}\label{r}
&\Big|\ss_{k\neq i,j}\int^t_0\p_{w_{i,+}}p_{jk}(w^{(m)}_+)(\p_sw^{(m)}_{i,+}+\la_j(w^{(m)}_+)\p_xw^{(m)}_{i,+})w^{(m)}_{k,+}(\psi^{(m)}_{t,j}(x,s)ds\Big|\no\\
=&\Big|\ss_{k\neq i,j}\int^t_0\p_{w_{i,+}}p_{jk}(w^{(m)}_+)(\p_sw^{(m)}_{i,+}-\p_sw^{(0)}_{i,+}+\la_j(w^{(m)}_+)(\p_xw^{(m)}_{i,+}-\p_xw^{(0)}_{i,+})\no\\
&+(\la_j(w^{(m)}_+)-\la_i(w^{(0)}_{+}))\p_xw^{(0)}_{i,+}(\psi^{(m)}_{t,j}(x,s),s))w^{(m)}_{k,+}(\psi^{(m)}_{t,j}(x,s),s)ds\Big|\no\\
\leqslant&M^7t^2+M^4t\int^t_0(s^3+(\psi^{(m)}_{t,j}(x,s)-\phi(s))^2)^{-\f{1}{6}}ds\no\\
\leqslant&M^8t(t^3+(x-\phi(t))^2)^{\f{2}{9}}.
\end{align}
As in \eqref{3-45}, for $|x-\phi(t)|\leqslant\ve^{\f{1}{2}}$, one has $|p_{jk}(w^{(m)}_{\pm}(x,t))|\leqslant M^2\ve^{\f{1}{3}}$.
This, together with \eqref{wjm+1e} and \eqref{r}, derives
\begin{align}\label{wj+}
|w^{(m+1)}_{j,+}(x,t)|&\leqslant M^4\ve^{\f{1}{3}}t(t^3+(x-\phi(t))^2)^{\f{1}{6}}+M^9t(t^3+(x-\phi(t))^2)^{\f{2}{9}}\no\\
&\leqslant t(t^3+(x-\phi(t))^2)^{\f{1}{6}},\quad\quad\quad j<i.
\end{align}
In the same way, we can arrive at
\begin{equation*}
|w^{(m+1)}_{j,-}(x,t)|\leqslant t(t^3+(x-\phi(t))^2)^{\f{1}{6}},\quad\quad\quad j>i.
\end{equation*}
Taking the first order derivative with respect to $x$ on both sides of \eqref{wjm+1t} yields
\begin{align*}
&\f{d}{ds}(\p_xw^{(m+1)}_{j,+}\circ\psi^{(m)}_{t,j}\p_x\psi^{(m)}_{t,j}(x,s)+\ss_{k\neq i,j}p_{jk}(w^{(m)}_+)\p_xw^{(m)}_{k,+}\circ\psi^{(m)}_{t,j}\p_x\psi^{(m)}_{t,j}(x,s))\no\\
=&\ss_{k\neq i,j}\f{d}{ds}p_{jk}(w^{(m)}_+)\p_xw^{(m)}_{k,+}\circ\psi^{(m)}_{t,j}\p_x\psi^{(m)}_{t,j}(x,s)
-\ss_{k\neq i,j}\p_xp_{jk}(w^{(m)}_+)\p_x\psi^{(m)}_{t,j}(x,s)\f{d}{ds}w^{(m)}_{k,+}\circ\psi^{(m)}_{t,j}.
\end{align*}
Due to $\p_xw^{(m+1)}_{k,+}(x,0)=\p_xw^{(m)}_{k,+}(x,0)=0$ for $k\neq i$ and \eqref{psi}, then
\begin{align}\label{wjm+}
\p_xw^{(m+1)}_{j,+}&(x,t)=-\ss_{k\neq i,j}p_{jk}(w^{(m)}_+(x,t))\p_xw^{(m)}_{k,+}(x,t)\no\\
&+\int^t_0\ss_{k\neq i,j}\ss_l\p_{w_{l,+}}p_{jk}(w^{(m)}_+)\f{d}{ds}w^{(m)}_{l,+}\circ\psi^{(m)}_{t,j}\p_xw^{(m)}_{k,+}\circ\psi^{(m)}_{t,j}\p_x\psi^{(m)}_{t,j}(x,s)ds\no\\
&-\int^t_0\ss_{k\neq i,j}\ss_l\p_{w_{l,+}}p_{jk}(w^{(m)}_+)\p_xw^{(m)}_{l,+}\circ\psi^{(m)}_{t,j}\p_x\psi^{(m)}_{t,j}(x,s)\f{d}{ds}w^{(m)}_{k,+}\circ\psi^{(m)}_{t,j}ds.
\end{align}
It follows from \eqref{wm-w0}-\eqref{ptwjm}, \eqref{pw0psin}, \eqref{h} and \eqref{suppsi} that
\begin{align}\label{+}
|\p_xw^{(m+1)}_{j,+}(x,t)|\leqslant& M^6\ve^{\f{1}{3}}(t^3+(x-\phi(t))^2)^{\f{1}{6}}+M^{10}(t^3+(x-\phi(t))^2)^{\f{2}{9}}\no\\
\leqslant& (t^3+(x-\phi(t))^2)^{\f{1}{6}},\quad j<i.
\end{align}
Combining with the equation in \eqref{wjm+1C}, this derives
\begin{equation*}
|\p_tw^{(m+1)}_{j,+}(x,t)|\leqslant M(t^3+(x-\phi(t))^2)^{\f{1}{6}},\quad j<i.
\end{equation*}
Similarly,
\begin{align*}
&|\p_xw^{(m+1)}_{j,-}(x,t)|\leqslant (t^3+(x-\phi(t))^2)^{\f{1}{6}},\quad j>i,\no\\
&|\p_tw^{(m+1)}_{j,-}(x,t)|\leqslant M(t^3+(x-\phi(t))^2)^{\f{1}{6}},\quad j>i.
\end{align*}

Next, we study the initial-boundary-value problem of $w^{(m+1)}_{j,-}$ for $j<i$,
\beq
\begin{cases}\label{wj-}
\p_tw^{(m+1)}_{j,-}+\la_j(w^{(m)}_-)\p_xw^{(m+1)}_{j,-}+\ss_{k\neq i,j}p_{jk}(w^{(m)}_-)(\p_tw^{(m)}_{k,-}+\la_j(w^{(m)}_-)\p_xw^{(m)}_{k,-})=0,\quad j<i,\\[2mm]
 w^{(m+1)}_{j,-}(x,t)|_{t=0}=0,
     \quad w^{(m+1)}_{j,-}(x,t)|_{x=\phi(t)}=w^{(m+1)}_{j,-}(\phi(t)_-,t),\quad j<i.
\end{cases}
\eeq
If the characteristics $\psi^{(m)}_{t,j}(x,s)$ intersects the time $\{t=0\}$, then as treated in the above we have
\begin{equation*}
|w^{(m+1)}_{j,-}(x,t)|\leqslant t(t^3+(x-\phi(t))^2)^{\f{1}{6}},\quad j<i.
\end{equation*}
If the characteristics $\psi^{(m)}_{t,j}(x,s)$ intersects the approximate shock curve at the point $(\phi(t_j),t_j)$ and $t_j>0$,
then the solution of \eqref{wj-} admits the following form
\begin{align}\label{wjme}
&w^{(m+1)}_{j,-}(x,t)\no\\
=&-\ss_{k\neq i,j}p_{jk}(w^{(m)}_-)w^{(m)}_{k,-}(x,t)+w^{(m+1)}_{j,-}(\phi(t_j)_-,t_j)+\ss_{k\neq i,j}p_{jk}(w^{(m)}_-)w^{(m)}_{k,-}(\phi(t_j)_-,t_j)\no\\
&+\ss_{k\neq i,j}\int^t_{t_j}\p_{w_{i,-}}p_{jk}(w^{(m)}_-)(\p_sw^{(m)}_{i,-}+\la_j(w^{(m)}_-)\p_xw^{(m)}_{i,-})w^{(m)}_{k,-}\circ\psi^{(m)}_{t,j}ds\no\\
&+\ss_{k\neq i,j}\ss_{l\neq i}\int^t_{t_j}\p_{w_{l,-}}p_{jk}(w^{(m)}_-)(\p_sw^{(m)}_{l,-}+\la_j(w^{(m)}_-)\p_xw^{(m)}_{l,-})w^{(m)}_{k,-}\circ\psi^{(m)}_{t,j}ds.
\end{align}
At first, we estimate $w^{(m+1)}_{j,-}(\phi(t_j)_-,t_j)$ for $j< i$.
It is derived from \eqref{w0j} and \eqref{wm-w0} that
\begin{align}\label{wij}
&|[w^{(m+1)}_i(t)]|=|w^{(m+1)}_{i,-}(\phi(t)_-,t)-w^{(m+1)}_{i,+}(\phi(t)_+,t)|\no\\
\leqslant&|w^{(m+1)}_{i,-}(\phi(t)_-,t)-w^{(0)}_{i,-}(\phi(t)_-,t)|+|w^{(0)}_{i,-}(\phi(t)_-,t)-w^{(0)}_{i,+}(\phi(t)_+,t)|\no\\
&+|w^{(m+1)}_{i,+}(\phi(t)_+,t)-w^{(0)}_{i,+}(\phi(t)_+,t)|\no\\
\leqslant&15t^{\f{1}{2}}.
\end{align}
Then combining with \eqref{RH} and \eqref{wij}, this yields
\begin{equation}\label{3-44}
|[w^{(m+1)}_j(t)]|\leqslant (15)^3Mt^{\f{3}{2}}.
\end{equation}
By the analogous estimates of \eqref{wjm+1e} and \eqref{wjm}, together with \eqref{3-44} and \eqref{wj+},
we obtain that for all $(x,t)\in\mathcal{B}_{\ve}$,
\begin{align*}
|w^{(m+1)}_{j,-}(\phi(t_j)_-,t_j)|&\leqslant|[w^{(m+1)}_{j}(t_j)]|+|w^{(m+1)}_{j,+}(\phi(t_j)_+,t_j)|\no\\
&\leqslant (15)^3Mt^{\f{3}{2}}+t(t^3+(\phi(t_j)_+-\phi(t_j))^2)^{\f{1}{6}}\no\\
&\leqslant \f{3}{4}M^2t(t^3+(x-\phi(t))^2)^{\f{1}{6}}.
\end{align*}
Thus, it follows from \eqref{r} and \eqref{wjme} that
\begin{equation*}
|w^{(m+1)}_{j,-}|\leqslant\f{4}{5}M^2t(t^3+(x-\phi(t))^2)^{\f{1}{6}}.
\end{equation*}

Taking the first order derivative with respect to $x$ on both sides of \eqref{wj-}, when
the characteristics $\psi^{(m)}_{t,j}(x,s)$ meets $\{t=0\}$, then as treated for \eqref{wjm+} we have
\begin{equation*}
|\p_xw^{(m+1)}_{j,-}(x,t)|\leqslant (t^3+(x-\phi(t))^2)^{\f{1}{6}},\quad j<i.
\end{equation*}
In addition, when the characteristics $\psi^{(m)}_{t,j}(x,s)$ intersects the approximate shock curve at the point $(\phi(t_j),t_j)$,
then
\begin{align}\label{pxwj}
&\p_xw^{(m+1)}_{j,-}(x,t)\no\\
=&-\ss_{k\neq i,j}(p_{jk}(w^{(m)}_-)\p_xw^{(m)}_{k,-})(x,t)+\p_xw^{(m+1)}_{j,-}(\phi(t_j)_-,t_j)+\ss_{k\neq i,j}(p_{jk}(w^{(m)}_-)\p_xw^{(m)}_{k,-})(\phi(t_j)_-,t_j)\no\\
&+\int^t_{t_j}\ss_{k\neq i,j}\ss_l\p_{w_{l,-}}p_{jk}(w^{(m)}_-\circ\psi^{(m)}_{t,j})\f{d}{ds}w^{(m)}_{l,-}\circ\psi^{(m)}_{t,j}\p_xw^{(m)}_{k,-}\circ\psi^{(m)}_{t,j}\p_x\psi^{(m)}_{t,j}ds\no\\
&-\int^t_{t_j}\ss_{k\neq i,j}\ss_l\p_{w_{l,-}}p_{jk}(w^{(m)}_-\circ\psi^{(m)}_{t,j})\p_xw^{(m)}_{l,-}\circ\psi^{(m)}_{t,j}\p_x\psi^{(m)}_{t,j}\f{d}{ds}w^{(m)}_{k,-}\circ\psi^{(m)}_{t,j}ds.
\end{align}
Next we estimate $\p_xw^{(m+1)}_{j,-}(\phi(t_j)_-,t_j)$ for $j<i$.
To this end, we treat $\f{d}{d t_j}w^{(m+1)}_{j,-}(\phi(t_j)_-,t_j)$ firstly.

Taking the first order derivative with respect to $t$ on both sides of \eqref{RH} yields
\begin{align}\label{3-46}
\f{d}{dt}[w^{(m+1)}_{j}(t)]=&(3\mathcal{F}_j[w^{(m+1)}_i(t)]^2+\p_{[w^{(m+1)}_i(t)]}\mathcal{F}_j[w^{(m+1)}_i(t)]^3)\f{d}{dt}[w^{(m+1)}_i(t)]\no\\
&+\p_{\langle w^{(m+1)}_i(t)\rangle}\mathcal{F}_j[w^{(m+1)}_i(t)]^3\f{d}{dt}\langle w^{(m+1)}_i(t)\rangle+O(t^{\f{3}{2}}),
\end{align}
where $\langle w^{(m+1)}_i(t)\rangle=\f{1}{2}(w^{(m+1)}_{i,-}(\phi(t)_-,t)+w^{(m+1)}_{i,+}(\phi(t)_+,t))$.
Note that
\begin{align*}
\Big|\f{d}{dt}[w^{(m+1)}_i(t)]\Big|&=\Big|\f{d}{dt}w^{(m+1)}_{i,-}(\phi(t)_-,t)-\f{d}{dt}w^{(m+1)}_{i,+}(\phi(t)_+,t)\Big|\no\\
&\leqslant\Big|\f{d}{dt}w^{(m+1)}_{i,-}-\f{d}{dt}w^{(0)}_{i,-}\Big|+\Big|\f{d}{dt}w^{(0)}_{i,-}
-\f{d}{dt}w^{(0)}_{i,+}\Big|+\Big|\f{d}{dt}w^{(m+1)}_{i,+}-\f{d}{dt}w^{(0)}_{i,+}\Big|
\end{align*}
and
\begin{align*}
\Big|\f{d}{dt}\langle w^{(m+1)}_i(t)\rangle\Big|&=\f{1}{2}\Big|\f{d}{dt}w^{(m+1)}_{i,-}(\phi(t)_-,t)+\f{d}{dt}w^{(m+1)}_{i,+}(\phi(t)_+,t)\Big|\no\\
&\leqslant\f{1}{2}\Big|\f{d}{dt}w^{(m+1)}_{i,-}-\f{d}{dt}w^{(0)}_{i,-}\Big|+\f{1}{2}\Big|\f{d}{dt}w^{(0)}_{i,-}
+\f{d}{dt}w^{(0)}_{i,+}\Big|+\f{1}{2}\Big|\f{d}{dt}w^{(m+1)}_{i,+}-\f{d}{dt}w^{(0)}_{i,+}\Big|.
\end{align*}
It follows from \eqref{jb}, \eqref{g0}, \eqref{g1}, \eqref{pxwm}, \eqref{ptwm}, \eqref{lam} and \eqref{wm+1-w0} that
\begin{align*}
\Big|\f{d}{dt}w^{(m+1)}_{i,+}-\f{d}{dt}w^{(0)}_{i,+}\Big|=&|\p_t(w^{(m+1)}_{i,+}-w^{(0)}_{i,+})+\phi'(t)\p_x(w^{(m+1)}_{i,+}-w^{(0)}_{i,+})|\no\\
\leqslant&|(\phi'(t)-\la_i(w^{(m)}_+))\p_x(w^{(m+1)}_{i,+}-w^{(0)}_{i,+})|+|(\la_i(w^{(m)}_+)-\la_i(w^{(0)}_+))\p_xw^{(0)}_{i,+}|\no\\
&+|\ss_{k\neq i}p_{ik}(w^{(m)}_+)(\p_tw^{(m)}_{k,+}+\la_i(w^{(m)}_{+})\p_xw^{(m)}_{k,+})|\no\\
\leqslant& M^7.
\end{align*}
Similarly, $\Big|\f{d}{dt}w^{(m+1)}_{i,-}-\f{d}{dt}w^{(0)}_{i,-}\Big|\leqslant M^7$ is obtained.
Thus, together with \eqref{wi0}, this derives
\beq\label{w}
\Big|\f{d}{d t}[w^{(m+1)}_i(t)]\Big|\leqslant Mt^{-\f{1}{2}},\quad\quad\Big|\f{d}{d t}\langle w^{(m+1)}_i(t)\rangle\Big|\leqslant M^8.
\eeq
Therefore, the estimate of \eqref{3-46} is
\begin{equation}\label{3-48}
  \Big|\f{d}{dt}[w^{(m+1)}_j(t)]\Big|\leqslant680M^2t^{\f{1}{2}}.
\end{equation}
By \eqref{RH}, \eqref{+}, \eqref{w}  and \eqref{3-48}, it holds that
\begin{align}\label{wj}
\Big|\f{d}{d t_j}w^{(m+1)}_{j,-}(\phi(t_j)_-,t_j)\Big|&=\Big|\f{d}{d t_j}[w^{(m+1)}_{j}(t_j)]+\f{d}{d t_j}w^{(m+1)}_{j,+}(\phi(t_j)_+,t_j)\Big|\no\\&\leqslant 800 M^2(t^3+(x-\phi(t))^2)^{\f{1}{6}}.
\end{align}
On the other hand, due to
\begin{equation*}
\begin{cases}
\p_{t_j} w^{(m+1)}_{j,-}(\phi(t_j)_-, t_j)+\phi'(t_j)\p_x w^{(m+1)}_{j,-}(\phi(t_j)_-, t_j)=\f{d}{d t_j}w^{(m+1)}_{j,-}(\phi(t_j)_-, t_j),\\[2mm]
\p_{t_j} w^{(m+1)}_{j,-}(\phi(t_j)_-, t_j)+\la_j(w^{(m)}_-(\phi(t_j)_-, t_j))\p_x w^{(m+1)}_{j,-}(\phi(t_j)_-, t_j)\\[2mm]
=-\ss_{k\neq i,j}(p_{jk}(w^{(m)}_-)(\p_{t_j} w^{(m)}_{k,-}+\la_j(w^{(m)}_-)\p_x w^{(m)}_{k,-}))(\phi(t_j)_-, t_j),
\end{cases}
\end{equation*}
then
\begin{align*}
&|\p_xw^{(m+1)}_{j,-}(\phi(t_j)_-,t_j)|\no\\
=&\Big|\f{1}{\phi'(t_j)-\la_j(w^{(m)}_-)}\Big(\f{d}{d t_j}w^{(m+1)}_{j,-}(\phi(t_j)_-,t_j)\no\\
&+\ss_{k\neq i,j}(p_{jk}(w^{(m)}_-)(\p_{t_j}w^{(m)}_{k,-}+\la_j\p_xw^{(m)}_{k,-})(\phi(t_j)_-,t_j))\Big)\Big|\no\\
\leqslant& 1000M^3(t^3+(x-\phi(t))^2)^{\f{1}{6}}+M^7\ve^{\f{1}{6}}(t^3+(x-\phi(t))^2)^{\f{1}{6}}\no\\
\leqslant&\f{1}{2}M^4(t^3+(x-\phi(t))^2)^{\f{1}{6}}.
\end{align*}
Returning to \eqref{pxwj}, by \eqref{h} and \eqref{suppsi}, we obtain
\begin{align}\label{pxwjm+1}
|\p_xw^{(m+1)}_{j,-}(x,t)|\leqslant&\f{1}{2}M^4(t^3+(x-\phi(t))^2)^{\f{1}{6}}+M^5\ve^{\f{1}{3}}(t^3+(x-\phi(t))^2)^{\f{1}{6}}\no\\
&+M^6(t^3+(x-\phi(t))^2)^{\f{2}{9}}\no\\
\leqslant& \f{3}{4}M^4(t^3+(x-\phi(t))^2)^{\f{1}{6}},\quad j<i.
\end{align}
Combining with the equation in \eqref{wj-} and \eqref{pxwjm+1}, one has
\begin{equation*}
|\p_t w^{(m+1)}_{j,-}(x,t)|\leqslant \f{3}{4}M^5(t^3+(x-\phi(t))^2)^{\f{1}{6}},\quad j<i.
\end{equation*}
Similarly,
\begin{align}
&|\p_xw^{(m+1)}_{j,+}(x,t)|\leqslant \f{3}{4}M^4(t^3+(x-\phi(t))^2)^{\f{1}{6}},\quad j>i,\no\\
&|\p_tw^{(m+1)}_{j,+}(x,t)|\leqslant \f{3}{4}M^5(t^3+(x-\phi(t))^2)^{\f{1}{6}},\quad j>i.\no
\end{align}
Thus we have finished the proof of Lemma \ref{diedai}.
$\hfill\space\square$

In order to obtain the convergence of the iteration scheme, we give the following contractive estimate.
\begin{lemma}
The map $\chi^{(m)}_{\ve}\longrightarrow\chi^{(m+1)}_{\ve}$ corresponding to
\begin{equation*}
(w^{(m)}_{1,\pm},\cdots,w^{(m)}_{i,\pm},\cdots,w^{(m)}_{n,\pm})\longmapsto(w^{(m+1)}_{1,\pm},\cdots,w^{(m+1)}_{i,\pm},\cdots,w^{(m+1)}_{n,\pm})
\end{equation*}
satisfies
\begin{align}\label{m}
&\max_{s\in[0,t]}(|w^{(m+1)}_{i,\pm}(\cdot,s)-w^{(m)}_{i,\pm}(\cdot,s)|+\ss_{j\neq i}|w^{(m+1)}_{j,\pm}(\cdot,s)-w^{(m)}_{j,\pm}(\cdot,s)|)\no\\
\leqslant&\f{1}{2}\max_{s\in[0,t]}(|w^{(m)}_{i,\pm}(\cdot,s)-w^{(m-1)}_{i,\pm}(\cdot,s)|+\ss_{j\neq i}|w^{(m)}_{j,\pm}(\cdot,s)-w^{(m-1)}_{j,\pm}(\cdot,s)|).
\end{align}
\end{lemma}
\pf
By the equations of the $(m+1)$-th iteration and the $m$-th iteration in the second equation of \eqref{im}, it holds that
\begin{align*}
    w^{(m+1)}_{i,\pm}&(x,t)-w^{(m)}_{i,\pm}(x,t)=-\int^t_0(\la_i(w^{(m)}_{\pm})-\la_i(w^{(m-1)}_{\pm}))\p_xw^{(0)}_{i,\pm}\circ\eta^{(m)}dt'\no\\
    &-\int^t_0(\la_i(w^{(m)}_{\pm})-\la_i(w^{(m-1)}_{\pm}))(\p_xw^{(m)}_{i,\pm}-\p_xw^{(0)}_{i,\pm})\circ\eta^{(m)}dt'\no\\
    &-\ss_{k\neq i}\int^t_0 (p_{ik}(w^{(m)}_{\pm})-p_{ik}(w^{(m-1)}_{\pm}))(\p_tw^{(m-1)}_{k,\pm}
    +\la_i(w^{(m-1)}_{\pm})\p_xw^{(m-1)}_{k,\pm})\circ\eta^{(m)}dt'\no\\
    &-\ss_{k\neq i}p_{ik}(w^{(m)}_{\pm}(x,t))(w^{(m)}_{k,\pm}(x,t)-w^{(m-1)}_{k,\pm}(x,t))\no\\
    &+\ss_{k\neq i}\ss_l\int^t_0(\p_{w_{l,\pm}}p_{ik}(w^{(m)}_{\pm})(\p_tw^{(m)}_{l,\pm}+\la_i(w^{(m)}_{\pm})\p_xw^{(m)}_{l,\pm}))\circ\eta^{(m)}(w^{(m)}_{k,\pm}-w^{(m-1)}_{k,\pm})dt'\no\\
    &-\ss_{k\neq i}\int^t_0(\la_i(w^{(m)}_{\pm})-\la_i(w^{(m-1)}_{\pm}))(p_{ik}(w^{(m)}_{\pm})\p_xw^{(m-1)}_{k,\pm})\circ\eta^{(m)}dt',
\end{align*}
where $w^{(m+1)}_{i,\pm}(x,0)-w^{(m)}_{i,\pm}(x,0)=0$ is used.
Thus, together with \eqref{wm-w0}-\eqref{ptwjm}, \eqref{w0etan}, \eqref{s} and \eqref{tt'}, we arrive at
\beq\label{n}
|w^{(m+1)}_{i,\pm}-w^{(m)}_{i,\pm}|\leqslant\Big(\f{12}{25}+C\ve^{\f{1}{6}}\Big)|w^{(m)}_{i,\pm}-w^{(m-1)}_{i,\pm}|+\ss_{j\neq i}C|w^{(m)}_{j,\pm}-w^{(m-1)}_{j,\pm}|,
\eeq
where and below $C=C(M)>0$ is a general constant independent of $m$.
Similarly, when $j<i$, based on \eqref{wm-w0}-\eqref{ptwjm}, \eqref{pw0psin} and \eqref{h},
one has that by first equation in \eqref{im}
\begin{align}\label{wm+1}
&|w^{(m+1)}_{j,+}(x,t)-w^{(m)}_{j,+}(x,t)|\no\\
=&\Big|-\int^t_0(\la_j(w^{(m)}_{+})-\la_j(w^{(m-1)}_{+}))\p_xw^{(m)}_{j,+}\circ\psi^{(m)}_{t,j}ds\no\\
&-\ss_{k\neq i,j}\int^t_0 (p_{jk}(w^{(m)}_{+})-p_{jk}(w^{(m-1)}_{+}))(\p_tw^{(m-1)}_{k,+}+\la_j(w^{(m-1)}_{+})\p_xw^{(m-1)}_{k,+})\circ\psi^{(m)}_{t,j}ds\no\\
&-\ss_{k\neq i,j}p_{jk}(w^{(m)}_{+}(x,t))(w^{(m)}_{k,+}(x,t)-w^{(m-1)}_{k,+}(x,t)) \no\\
    &+\ss_{k\neq i,j}\ss_l\int^t_0(\p_{w_{l,+}}p_{jk}(w^{(m)}_{+})(\p_tw^{(m)}_{l,+}+\la_j(w^{(m)}_{+})\p_xw^{(m)}_{l,+}))\circ\psi^{(m)}_{t,j}(w^{(m)}_{k,+}-w^{(m-1)}_{k,+})ds\no\\
    &-\ss_{k\neq i,j}\int^t_0(\la_j(w^{(m)}_{+})-\la_j(w^{(m-1)}_{+}))(p_{jk}(w^{(m)}_{+})\p_xw^{(m-1)}_{k,+})\circ\psi^{(m)}_{t,j}ds\Big|\no\\
    \leqslant&C\ve^{\f{1}{6}}|w^{(m)}_{i,+}-w^{(m-1)}_{i,+}|+\ss_{k\neq i,j}C\ve^{\f{1}{6}}|w^{(m)}_{k,+}-w^{(m-1)}_{k,+}|+C\ve^{\f{1}{6}}|w^{(m)}_{j,+}-w^{(m-1)}_{j,+}|.
\end{align}
For the first equation in \eqref{ibm}, due to
\begin{align}\label{p}
&|(w^{(m+1)}_{j,-}-w^{(m)}_{j,-})(\phi(t_j)_-,t_j)|
=|[w^{(m+1)}_j]+w^{(m+1)}_{j,+}(\phi(t_j)_+,t_j)-[w^{(m)}_j]-w^{(m)}_{j,+}(\phi(t_j)_+,t_j)|\no\\
&\quad \leqslant |\mathcal{F}_j[w^{(m+1)}_i]^3-\mathcal{F}_j[w^{(m)}_i]^3|+|w^{(m+1)}_{j,+}-w^{(m)}_{j,+}|\no\\
&\quad \leqslant C\ve|w^{(m+1)}_{i,\pm}-w^{(m)}_{i,\pm}|+|w^{(m+1)}_{j,+}-w^{(m)}_{j,+}|,
\end{align}
then together with \eqref{n}, \eqref{wm+1} and \eqref{p}, we obtain \eqref{m}.
\ef
Set $z=x-\phi(t)$ and
\beq\label{wz}
\mathrm{w}_i(z,t)=w_i(x,t),\quad\quad\mathrm{w}_j(z,t)=w_j(x,t),\quad j\neq i.
\eeq
Then the related domain becomes $\mathrm{B}_{\ve}=([-M\ve,M\ve]\backslash\{0\})\times[0,\ve]$.

By the compression mapping principle, there exists a vector function $(\mathrm{w}_1,\cdots,\mathrm{w}_i,\cdots,\mathrm{w}_n)\in C(\mathrm{B}_{\ve})
\cap C^1(([-M\ve,M\ve]\backslash\{0\})\times (0,\ve))$ such that $(\mathrm{w}^{(m)}_1,\cdots,\mathrm{w}^{(m)}_i,\cdots,\mathrm{w}^{(m)}_n)\longrightarrow(\mathrm{w}_1,\cdots,\mathrm{w}_i,\cdots,\mathrm{w}_n)$ uniformly in $\mathrm{B}_{\ve}$. Due to the boundedness of first order derivatives,
according to the Banach-Alaoglu theorem and the lower semi-continuity, we obtain
\begin{align*}
(t^3+z^2)^{\f{1}{6}}\mathrm{w}^{(m)}_i&\rightharpoonup(t^3+z^2)^{\f{1}{6}}\mathrm{w}_i,\quad\text{in}\quad W^{1,\infty}(([-M\ve,M\ve]\backslash\{0\})\times (0,\ve))\quad\text{weak-*},\no\\
\mathrm{w}^{(m)}_j&\rightharpoonup\mathrm{w}_j,~\quad j\neq i,~~\quad\text{in}\quad W^{1,\infty}(([-M\ve,M\ve]\backslash\{0\})\times (0,\ve))\quad\text{weak-*}.
\end{align*}
Thus, $(\mathrm{w}_1,\cdots,\mathrm{w}_i,\cdots,\mathrm{w}_n)$ constitutes the solution of \eqref{f} by \eqref{jb}.

\section{The existence and uniqueness of the shock solution}\label{Sec-4}
For the existence and uniqueness of the solution corresponding to \eqref{f}, \eqref{jb} has been assumed.
In this section, we will close this assumption and determine the real shock curve $x=\vp(t)$.

From now on, we consider the slope $\si(t)=\vp'(t)$ of the shock curve $x=\vp(t)$. Due to the Corollary 17.3
in \cite{Smoller} and \eqref{RH-0},
we can write
\begin{equation*}
\si(t)=\la_i(w_+)+\f{1}{2}\p_{w_{i,+}}\la_i(w_+)[w_i]+\ss_{j\neq i}\mathrm{F_j}(w_+)[w_j]+\sum^n_{i,j=1}\mathrm{F_{ij}}(w_{\pm})[w_i][w_j],
\end{equation*}
where $\mathrm{F_j}(w_+)$ and $\mathrm{F_{ij}}(w_{\pm})$ are smooth functions, and
\begin{align*}
\la_i(w_+)=&w_{i,+}+\ss_{j\neq i}\p_{w_{j,+}}\la_i(0)w_{j,+}
+O((w_{i,+})^2+\ss_{j\neq i}(w_{j,+})^2).
\end{align*}
In addition, $w_{i,+}=\langle w_i\rangle-\f{1}{2}[w_i]$ holds.
Thus, together with $\p_{w_i}\la_i(0)=1$, $\si(t)$ has the following expression
\begin{align}\label{sit}
  \si (t) =&\langle w_i\rangle+\ss_{j\neq i}\p_{w_{j,+}}\la_i(0)w_{j,+}
  +\ss_{j\neq i}\mathrm{F_j}(0)[w_j]+\f{1}{2}\p^2_{w_{i,+}}\la_i(0)(\langle w_i \rangle -\f{1}{2}[w_i])[w_i]\no\\
  &+\sum^n_{i,j=1}\mathrm{F_{ij}}(0)[w_i][w_j]+O((\langle w_i\rangle-\f{1}{2}[w_i])^2+\ss_{j\neq i}(w_{j,+})^2).
\end{align}
By \eqref{RH}, \eqref{w0j},  \eqref{wi0m}, \eqref{wm-w0} and \eqref{wjm}, we have
\begin{equation}\label{Sun-3.1}
|\si(t)|\leqslant M^4 t.
\end{equation}
Integrating \eqref{Sun-3.1} over $[0,t]$ yields the first estimate in \eqref{jb}.
Besides, taking the first order derivative with respect to $t$ in \eqref{sit}, by \eqref{pxwjm}, \eqref{ptwjm}, \eqref{w} and \eqref{wj},
we have
\begin{align*}
  |\si'(t)|
  =&\Big|\f{d}{dt}\langle w_i\rangle+\ss_{j\neq i}\p_{w_{j,+}}\la_i(0)\f{d}{dt}w_{j,+}+\ss_{j\neq i}\mathrm{F_j}(0)\f{d}{dt}[w_j]\no\\
  &+\f{1}{2}\p^2_{w_{i,+}}\la_i(0)\Big(\f{d}{dt}\langle w_i \rangle -\f{1}{2}\f{d}{dt}[w_i]\Big)[w_i]\no\\
  &+\f{1}{2}\p^2_{w_{i,+}}\la_i(0)(\langle w_i \rangle -\f{1}{2}[w_i])\f{d}{dt}[w_i]\no\\
  &+\sum^n_{i,j=1}\mathrm{F_{ij}}(0)\Big(\f{d}{dt}[w_i][w_j]+[w_i]\f{d}{dt}[w_j]\Big)\no\\
  &+O((\langle w_i\rangle-\f{1}{2}[w_i])\Big(\f{d}{dt}\langle w_i\rangle-\f{1}{2}\f{d}{dt}[w_i]\Big)+w_{j,+}\f{d}{dt}w_{j,+})\Big|\no\\
\leqslant& \f{1}{2}M^9.
\end{align*}
Thus we have finished the closure of the bootstrap assumptions in \eqref{jb}.

Denote $\{\varphi^{(m)}(t)\}_{m\geqslant0}$ by the iterative approximate shock curve, which is uniformly bounded in $W^{2,\infty}(0,\ve)$
and satisfies \eqref{jb}. Then we have
\begin{lemma}\label{3.8}
There exists a unique uniform limitation $\varphi(t)$ such that $\varphi(t)$ satisfies \eqref{jb}
corresponding to the system \eqref{f}, with the initial data introduced in \eqref{2-19}. Besides, the uniqueness of the entropy solution $w$ for the system \eqref{f} holds.
\end{lemma}
\pf
Firstly, we define a sequence of $\varphi^{(m)}$ for $m\geqslant0$ as
\begin{equation}\label{isi}
  \varphi^{(m+1)}(t)=\int^t_0\si^{(m)}(\tau)d\tau,\quad\quad\varphi^{(0)}(t)=0.
\end{equation}
By adjusting the domain as in \eqref{wz}, from \eqref{RH}, \eqref{w0j}, \eqref{wi0m} and \eqref{sit}, we have
\begin{align}\label{dotsii+1}
|\si^{(m+1)}(t)-\si^{(m)}(t)|\leqslant&(1+C\ve^{\f{1}{2}})|\langle \mathrm{w}^{(m+1)}_i\rangle-\langle \mathrm{w}^{(m)}_i\rangle|+C\ve^{\f{1}{2}}|[\mathrm{w}^{(m+1)}_i]-[\mathrm{w}^{(m)}_i]|\no\\
&+\ss_{j\neq i}C|\mathrm{w}^{(m+1)}_{j,\pm}-\mathrm{w}^{(m)}_{j,\pm}|.
\end{align}
Since $\varphi^{(m)}(t)$ satisfies \eqref{jb}, one can know that $(\mathrm{w}_1,\cdots,\mathrm{w}_i,\cdots,\mathrm{w}_n)^{(m)}$
constitutes the solution of \eqref{f} from the discussion of Section \ref{Sec-3}.
By using \eqref{sit} and \eqref{isi}, we also obtain $\varphi^{(m+1)}(t)$ satisfies \eqref{jb}. This means
that the map $\varphi^{(m)}(t)\longrightarrow\varphi^{(m+1)}(t)$ is stable. If we can prove $\si^{(m)}\longrightarrow\si(t)$ as $\varphi^{(m)}(t)\longrightarrow\varphi(t)$, then the limit function $\si(t)$ solves \eqref{RH-0},
the real shock curve $x=\varphi(t)$ will be constructed.

Note that for $j\neq i$ and $(z,t)\in\mathrm{B}_{\ve}$, we introduce
\begin{equation}\label{4-5}
(\mathrm{w}^{(m)}_1,\cdots,\mathrm{w}^{(m)}_i,\cdots,\mathrm{w}^{(m)}_n)(z,t)=(w^{(m)}_1,\cdots,w^{(m)}_i,\cdots,w^{(m)}_n)(z+\varphi^{(m)}(t),t).
\end{equation}
Thus, \eqref{RH} are satisfied for each $m\geqslant 0$ and $t\in[0,\ve]$. Furthermore, $\si^{(m)}(t)$ can be expressed by the
above variables in the left hand side of \eqref{4-5}. Then it follows from \eqref{f} in the coordinate $(z,t)$ of \eqref{4-5}
that
\begin{equation*}
\left\{
\begin{array}{ll}
(\p_t+(\la^{(m)}_i-\si^{(m)})\p_z)\mathrm{w}^{(m)}_{i,\pm}+\ss_{k\neq i}p_{ik}(\mathrm{w}^{(m)}_{\pm})(\p_t+(\la^{(m)}_i-\si^{(m)})\p_z)\mathrm{w}^{(m)}_{k,\pm}=0,\\[2mm]
(\p_t+(\la^{(m)}_j-\si^{(m)})\p_z)\mathrm{w}^{(m)}_{j,\pm}+\ss_{k\neq i,j}p_{jk}(\mathrm{w}^{(m)}_{\pm})(\p_t+(\la^{(m)}_j-\si^{(m)})\p_z)\mathrm{w}^{(m)}_{k,\pm}=0.
\end{array}
\right.
\end{equation*}
\\
In addition, the equation of $\mathrm{w}^{(m+1)}_{i,\pm}-\mathrm{w}^{(m)}_{i,\pm}$
can be written as
\begin{align}\label{4-3}
&(\p_t+(\la^{(m+1)}_i-\si^{(m+1)})\p_z)(\mathrm{w}^{(m+1)}_{i,\pm}-\mathrm{w}^{(m)}_{i,\pm}+\ss_{k\neq i}p_{ik}(\mathrm{w}^{(m+1)}_{\pm})(\mathrm{w}^{(m+1)}_{k,\pm}-\mathrm{w}^{(m)}_{k,\pm}))\no\\
=&-(\la^{(m+1)}_i-\si^{(m+1)}-\la^{(m)}_i+\si^{(m)})(\p_z\mathrm{w}^{(m)}_{i,\pm}+\ss_{k\neq i}p_{ik}(\mathrm{w}^{(m+1)}_{\pm})\p_zw^{(m)}_{k,\pm})\no\\
&+\ss_{k\neq i}\ss_l\p_{\mathrm{w}_{l,\pm}}p_{ik}(\mathrm{w}^{(m+1)}_{\pm})(\p_t\mathrm{w}^{(m+1)}_{l,\pm}
+(\la^{(m+1)}_i-\si^{(m+1)})\p_z\mathrm{w}^{(m+1)}_{l,\pm})(\mathrm{w}^{(m+1)}_{k,\pm}-\mathrm{w}^{(m)}_{k,\pm})\no\\
&-\ss_{k\neq i}\ss_l\p_{\mathrm{w}_{l,\pm}}p_{ik}(\mathrm{w}^{(m)}_{\pm})(\mathrm{w}^{(m+1)}_{l,\pm}
-\mathrm{w}^{(m)}_{l,\pm})(\p_t\mathrm{w}^{(m)}_{k,\pm}+(\la^{(m)}_i-\si^{(m)})\p_z\mathrm{w}^{(m)}_{k,\pm}).
\end{align}
Applying the characteristics method, in terms of \eqref{wm-w0}-\eqref{ptwjm}, \eqref{w0etan}
and \eqref{wz}, one has that by $\mathrm{w}^{(m+1)}_{i,\pm}(x,0)-\mathrm{w}^{(m)}_{i,\pm}(x,0)=0$ and $\mathrm{w}^{(m+1)}_{j,\pm}(x,0)-\mathrm{w}^{(m)}_{j,\pm}(x,0)=0$ for $j\neq i$,
\begin{equation*}
\Big|\int^t_0\p_z\mathrm{w}^{(m)}_{i,\pm}\circ(\eta^{(m+1)}-\varphi^{(m+1)})dt'\Big|
=\Big|\int^t_0\p_xw^{(m)}_{i,\pm}\circ\eta^{(m+1)}dt'\Big|\leqslant\f{24}{55}+C\ve^{\f{1}{6}}.
\end{equation*}
Through solving \eqref{4-3}, we arrive at
\begin{align*}
|\mathrm{w}^{(m+1)}_{i,\pm}-\mathrm{w}^{(m)}_{i,\pm}|\leqslant&\Big(\f{12}{25}+C\ve^{\f{1}{8}}\Big)|\mathrm{w}^{(m+1)}_{i,\pm}
-\mathrm{w}^{(m)}_{i,\pm}|+\Big(\f{24}{55}+C\ve^{\f{1}{8}}\Big)|\si^{(m+1)}-\si^{(m)}|\no\\
&+\ss_{j\neq i}\Big(\f{24}{55}C+C\ve^{\f{1}{8}}\Big)|\mathrm{w}^{(m+1)}_{j,\pm}-\mathrm{w}^{(m)}_{j,\pm}|.
\end{align*}
Therefore,
\beq\label{mw}
|\mathrm{w}^{(m+1)}_{i,\pm}-\mathrm{w}^{(m)}_{i,\pm}|\leqslant\f{17}{20}|\si^{(m+1)}-\si^{(m)}|+\ss_{j\neq i}\f{17}{20}C|\mathrm{w}^{(m+1)}_{j,\pm}-\mathrm{w}^{(m)}_{j,\pm}|.
\eeq
Similarly, the equation of $\mathrm{w}^{(m+1)}_{j,\pm}-\mathrm{w}^{(m)}_{j,\pm}$ for $j\neq i$ can be described as
\begin{align*}
&(\p_t+(\la^{(m+1)}_j-\si^{(m+1)})\p_z)(\mathrm{w}^{(m+1)}_{j,\pm}-\mathrm{w}^{(m)}_{j,\pm}+\ss_{k\neq i,j}p_{jk}(\mathrm{w}^{(m+1)}_{\pm})(\mathrm{w}^{(m+1)}_{k,\pm}-\mathrm{w}^{(m)}_{k,\pm}))\no\\
=&-(\la^{(m+1)}_j-\si^{(m+1)}-\la^{(m)}_j+\si^{(m)})(\p_z\mathrm{w}^{(m)}_{j,\pm}+\ss_{k\neq i,j}p_{jk}(\mathrm{w}^{(m+1)}_{\pm})\p_zw^{(m)}_{k,\pm})\no\\
&+\ss_{k\neq i,j}\ss_l\p_{\mathrm{w}_{l,\pm}}p_{jk}(\mathrm{w}^{(m+1)}_{\pm})(\p_t\mathrm{w}^{(m+1)}_{l,\pm}+(\la^{(m+1)}_j-\si^{(m+1)})\p_z\mathrm{w}^{(m+1)}_{l,\pm})(\mathrm{w}^{(m+1)}_{k,\pm}-\mathrm{w}^{(m)}_{k,\pm})\no\\
&-\ss_{k\neq i,j}\ss_l\p_{\mathrm{w}_{l,\pm}}p_{jk}(\mathrm{w}^{(m)}_{\pm})(\mathrm{w}^{(m+1)}_{l,\pm}-\mathrm{w}^{(m)}_{l,\pm})(\p_t\mathrm{w}^{(m)}_{k,\pm}+(\la^{(m)}_j-\si^{(m)})\p_z\mathrm{w}^{(m)}_{k,\pm}),
\end{align*}
and then it follows from \eqref{pw0psin} that
\begin{align}\label{mwj}
|\mathrm{w}^{(m+1)}_{j,\pm}-\mathrm{w}^{(m)}_{j,\pm}|\leqslant&C\ve^{\f{4}{9}}|\mathrm{w}^{(m+1)}_{i,\pm}-\mathrm{w}^{(m)}_{i,\pm}|+C\ve^{\f{4}{9}}|\si^{(m+1)}-\si^{(m)}|\no\\
&+\ss_{k\neq i,j}C\ve^{\f{1}{6}}|\mathrm{w}^{(m+1)}_{k,\pm}-\mathrm{w}^{(m)}_{k,\pm}|+C\ve^{\f{4}{9}}|\mathrm{w}^{(m+1)}_{j,\pm}-\mathrm{w}^{(m)}_{j,\pm}|.
\end{align}
Thus, combining with \eqref{mw} and \eqref{mwj} yields
\beq\label{mwi}
|\mathrm{w}^{(m+1)}_{i,\pm}-\mathrm{w}^{(m)}_{i,\pm}|+\ss_{j\neq i}\ve^{-\f{1}{12}}|\mathrm{w}^{(m+1)}_{j,\pm}-\mathrm{w}^{(m)}_{j,\pm}|\leqslant\f{15}{16}|\si^{(m+1)}-\si^{(m)}|.
\eeq
In terms of \eqref{isi}, \eqref{dotsii+1} and \eqref{mwi}, we have
\begin{align}\label{dotsii+1-}
  &\sup_{t\in[0,\ve]}|\si^{(m+2)}(t)-\si^{(m+1)}(t)|\no\\
   \leqslant&\sup_{t\in[0,\ve]}(1+C\ve^{\f{1}{2}})|\mathrm{w}^{(m+1)}_{i,\pm}- \mathrm{w}^{(m)}_{i,\pm}|+\ss_{j\neq i}C|\mathrm{w}^{(m+1)}_{j,\pm}-\mathrm{w}^{(m)}_{j,\pm}|\no\\
   \leqslant&\sup_{t\in[0,\ve]}(1+C\ve^{\f{1}{2}}+C\ve^{\f{1}{12}})(|\mathrm{w}^{(m+1)}_{i,\pm}-\mathrm{w}^{(m)}_{i,\pm}|+\ss_{j\neq i}\ve^{-\f{1}{12}}|\mathrm{w}^{(m+1)}_{j,\pm}-\mathrm{w}^{(m)}_{j,\pm}|)\no\\
   \leqslant&\f{19}{20}\sup_{t\in[0,\ve]}|\si^{(m+1)}(t)-\si^{(m)}(t)|.
\end{align}
Thus, there exists a unique limitation $\varphi(t)$ in $\mathrm{W}^{1,\infty}(0,\ve)$ and inherits the bound \eqref{jb}. In addition, $\si^{(m)}\longrightarrow\si(t)\in C^0(0,\ve)$ as $m\longrightarrow\infty$, furthermore, $\si(t)=\varphi'(t)$ is the solution of \eqref{RH-0}.

Since the above $\varphi(t)$ satisfies \eqref{jb}, the uniqueness of the entropy solution of \eqref{f} has been proved.
If there are two different shock curves $\varphi^{(1)}$ and $\varphi^{(2)}$ satisfying \eqref{jb}.
Then $(w_1,\cdots,w_i,\cdots,w_n)^{(m)}(m=1,2)$ are $C^1_{x,t}(\mathcal{B}_{\ve})$ smooth solutions of \eqref{f}
with the initial data \eqref{2-19}.
Thus, it follows from \eqref{mwi} and \eqref{dotsii+1-} that
\begin{equation*}
|\mathrm{w}^{(2)}_{i,\pm}-\mathrm{w}^{(1)}_{i,\pm}|+\ss_{j\neq i}\ve^{-\f{1}{12}}|\mathrm{w}^{(2)}_{j,\pm}-\mathrm{w}^{(1)}_{j,\pm}|=|\si^{(2)}-\si^{(1)}|=0,\quad t\in[0,\ve].
\end{equation*}
Furthermore, by $\varphi^{(m)}(0)=0 (m=1,2)$, $\varphi^{(1)}=\varphi^{(2)}$ holds.
In conclusion, we have finished the proof of Lemma \ref{3.8}.
\ef

\section{Proof of the main theorem}\label{Sec-5}
\noindent{\bf  Proof of Theorem \ref{main1}}\quad
By Section \ref{Sec-2}, we have obtained that the $i$-simple wave solution of \eqref{isimple} has $C^{\f{1}{3}}$ H\"older
regularity near the blowup point $(0,0)$. In addition, under the condition \eqref{jb},
by the characteristics method, the estimates \eqref{wm-w0}-\eqref{ptwjm} for \eqref{f} are obtained
in Section \ref{Sec-3}. On the other hand,  we have constructed a unique $C^2$ regular
shock curve $x=\vp(t)$ in Lemma \ref{3.8}. Meanwhile, the weak entropy solution $w$ of \eqref{f} satisfying \eqref{1.12}
is derived. Naturally, \eqref{1.12-0} is proved.
\ef

\vskip 0.3 true cm

{\bf \color{blue}{Conflict of Interest Statement:}}

\vskip 0.2 true cm

{\bf The authors declare that there is no conflict of interest in relation to this article.}

\vskip 0.2 true cm
{\bf \color{blue}{Data availability statement:}}

\vskip 0.2 true cm

{\bf  Data sharing is not applicable to this article as no data sets are generated
during the current study.}

\vskip 0.3 true cm

\end{document}